\documentclass[nameyear,seceqn,dvips,aos]{arxims}
\usepackage{dcolumn,mathbh}
\usepackage{graphicx}
\usepackage{hyperref}

\volume{32}
\issue{2}
\pubyear{2004}
\firstpage{407}
\lastpage{451}
\doi{10.1214/009053604000000067}

\makeatletter

\newcommand{\Sig}{\Sigma}

\newcommand{\lam}{\lambda}

\newcommand{\be}{\begin{eqnarray*}}
\newcommand{\ee}{\end{eqnarray*}}
\def\bx {{\bf x}}
\def\by {{\bf y}}
\def\br {{\bf r}}
\newcommand{\G}{{\mathcal G}}

\newproclaim{rem}{Remark}

\newcommand{\bfeta}{{\bolds{\eta}}}
\newcommand{\bdelta}{\bolds{\delta}}
\newcommand{\bbeta}{\bolds{\beta}}

\newcommand{\bmu}{{\bolds{\mu}}}
\newcommand{\bepsi}{\bolds{\varepsilon}}

\newtheorem{lemma}{Lemma}

\newproclaim{las}{Lasso modification}
\newproclaim{pos}{Positive Lasso}
\newproclaim{con}{Constraint}
\newproclaim{conI}{Constraint I}
\newproclaim{conII}{Constraint II}
\newproclaim{conIII}{Constraint III}
\newproclaim{conv}{Conventions}
\newproclaim{stag}{Stagewise modification}
\newproclaim{posit}{Positive cone condition}

\newcommand{\sign}{\operatorname{sign}}
\newcommand{\pe}{\operatorname{pe}}
\newcommand{\cov}{\operatorname{cov}}
\newcommand{\Cos}{\operatorname{Cos}}
\newcommand{\Sin}{\operatorname{Sin}}
\newcommand{\Proj}{\operatorname{Proj}}
\makeatother

\begin{document}
\begin{frontmatter}

\title{Least Angle Regression}
\runtitle{Least Angle Regression}

\begin{aug}
\author[A]{\fnms{Bradley} 
\snm{Efron}\corref{}\protect\thanksref{T1}\ead[label=e1]{brad@stat.stanford.edu}},
\thankstext[1]{T1}{Supported in part by NSF Grant DMS-00-72360 and
NIH Grant 8R01-EB002784.}
\author[B]{\fnms{Trevor} \snm{Hastie}\protect\thanksref{T2}},
\thankstext[2]{T2}{Supported in part by NSF Grant DMS-02-04162 and
NIH  Grant  R01-EB0011988-08.}
\author[C]{\fnms{Iain} 
\snm{Johnstone}\protect\thanksref{T3}}
\thankstext[3]{T3}{Supported in part by NSF Grant DMS-00-72661
and NIH Grant R01-EB001988-08.}
and
\author[D]{\fnms{Robert} 
\snm{Tibshirani}\protect\thanksref{T4}}
\thankstext[4]{T4}{Supported in part by NSF Grant DMS-99-71405
and  NIH Grant 2R01-CA72028.}
\runauthor{Efron,  Hastie,  Johnstone and  Tibshirani}
\affiliation{Stanford University}
\address[A]{Department of Statistics\\
Stanford University\\
Sequoia Hall\\
Stanford, California 94305-4065\\
USA\\
\printead{e1}} 
\end{aug}

\received{\smonth{3} \syear{2002}}
\revised{\smonth{1} \syear{2003}}

\begin{abstract}
The purpose of model selection algorithms such as {\it All Subsets},
{\em Forward Selection} and {\em Backward Elimination} is to choose
a linear model on the basis of the same set of data to which the
model will be applied. Typically we have available a large
collection of possible covariates from which we hope to select a
  parsimonious set for the efficient prediction of a response
  variable. {\em Least Angle Regression} (LARS), a new model
  selection algorithm, is a useful and less greedy version of
  traditional forward selection methods. Three main properties are
  derived: (1) A simple modification of the LARS algorithm implements
  the Lasso, an attractive version of ordinary least squares that
  constrains the sum of the absolute regression coefficients; the LARS
  modification calculates all possible Lasso estimates for a given
  problem, using an order of magnitude less computer time than previous
  methods. (2) A different LARS modification efficiently implements
  Forward Stagewise linear regression, another promising new model
  selection method; this connection explains the similar numerical
  results previously observed for the Lasso and Stagewise, and helps us
  understand the properties of both methods, which are seen as
  constrained versions of the simpler LARS algorithm. (3) A simple
  approximation for the degrees of freedom of a LARS estimate is
  available, from which we derive a $C_p$ estimate of prediction error;
  this allows a principled choice among the range of possible LARS
  estimates. LARS and its variants are computationally efficient: the paper 
describes a~publicly available algorithm that requires only the same order of
  magnitude of computational effort as ordinary least squares applied
  to the full set of covariates.
\end{abstract}

\setattribute{keyword}{AMS}{AMS 2000 subject classification.}

                      
\end{frontmatter}


\section{Introduction.}\label{sec:sec1} 
Automatic model-building algorithms
are familiar, and sometimes notorious, in the linear model literature:
Forward Selection, Backward Elimination, All Subsets regression and
various combinations are used to automatically produce ``good'' linear
models for predicting a response $y$ on the basis of some measured
covariates $x_1, x_2, \ldots, x_m$. Goodness is often defined in terms
of prediction accuracy, but parsimony is another important criterion:
simpler models are preferred for the sake of scientific insight into
the $x-y$ relationship. Two promising recent model-building
algorithms,\ the Lasso and Forward Stagewise linear regression, will be
discussed here, and motivated in terms of a computationally simpler
method called Least Angle Regression.

Least Angle Regression (LARS) relates to the classic
model-selection method known as Forward Selection, or ``forward
stepwise regression,'' described in \citeauthor{W80} [(\citeyear{W80}),
Section~8.5]:
given a collection of possible predictors, we select the one having
largest absolute correlation with the response $y$, say $x_{j_1}$, and
perform simple linear regression of $y$ on $x_{j_1}$. This leaves a
residual vector orthogonal to $x_{j_1}$, now considered to be the
response. We project the other predictors orthogonally to $x_{j_1}$
and repeat the selection process. After $k$ steps this results in a set
of predictors $x_{j_1}, x_{j_2},\ldots, x_{j_k}$ that are then used
in the usual way to construct a $k$-parameter linear model. Forward
Selection is an aggressive fitting technique that can be overly
greedy, perhaps eliminating at the second step useful predictors that
happen to be correlated with $x_{j_1}$.

Forward Stagewise, as described below, is a much more cautious version
of Forward Selection, which may take thousands of tiny steps as it
moves toward a final model. It turns out, and this was the original
motivation for the LARS algorithm, that a simple formula allows
Forward Stagewise to be implemented using fairly large steps, though
not as large as a classic Forward Selection, greatly reducing the
computational burden.  The geometry of the algorithm, described in
Section 2, suggests the name ``Least Angle Regression.'' 
It then
happens that this same geometry applies to another, seemingly quite
different, selection method called the Lasso [\citet{Ti96}]. The
LARS--Lasso--Stagewise connection is conceptually as well as
computationally useful. The Lasso is described next, in terms of the
main example used in this paper.

Table~\ref{tab1}   shows a small part of the data for our main example.

\begin{table}[b]
\tabcolsep=6.5pt
\caption{Diabetes study: $442$ diabetes patients
were measured on $10$ baseline variables; a prediction model was
desired for the response variable, a measure of disease progression
one year after baseline}
\label{tab1}
\begin{tabular}{@{}cccccccccccc@{}}    \hline
 & {\bf AGE} & {\bf SEX} & {\bf BMI} &  {\bf BP} &
\multicolumn{6}{@{}c@{}}{{\bf Serum measurements}} & {\bf Response}\\
\cline{6-11}
\noalign{}
{\bf Patient}       & ${\bf x_1}$ & ${\bf x_2}$ & ${\bf x_3}$ &  
${\bf x_4}$ & \multicolumn{1}{c}{${\bf x_5}$} & ${\bf x_6}$ &  
${\bf x_7}$ &  
${\bf x_8}$ &  
${\bf x_9}$ & \multicolumn{1}{c}{${\bf x_{10}}$} &   {\bf y} \\
\hline
\phantom{00}1&59&2&32.1&101&157&~93.2&38&4&4.9&~87&151\\ 
\phantom{00}2&48&1&21.6&~87&183&103.2&70&3&3.9&~69&~75\\ 
\phantom{00}3&72&2&30.5&~93&156&~93.6&41&4&4.7&~85&141\\ 
\phantom{00}4&24&1&25.3&~84&198&131.4&40&5&4.9&~89&206\\ 
\phantom{00}5&50&1&23.0&101&192&125.4&52&4&4.3&~80&135\\ 
\phantom{00}6&23&1&22.6&~89&139&~64.8&61&2&4.2&~68&~97\\ 
\phantom{00}    \vdots &      \vdots &  \vdots &   \vdots &    \vdots &   \vdots &   \vdots &    \vdots &     \vdots &    \vdots &  \vdots &       \vdots \\
441&36&1&30.0&~95&201&125.2&42&5&5.1&~85&220\\ 
442&36&1&19.6&~71&250&133.2&97&3&4.6&~92&~57\\ 
\hline
\end{tabular}
\end{table}

Ten baseline variables, age, sex, body mass index, average
blood pressure and six blood serum measurements, 
were obtained for
each of $n = 442$ diabetes patients, as well as the response of
interest, a quantitative measure of disease progression one year after
baseline. The statisticians were asked to construct a model that
predicted response $y$ from covariates $x_1, x_2, \ldots,
x_{10}$. Two hopes were evident here, that the model would produce
accurate baseline predictions of response for future patients and
that the form of the model would suggest which covariates were
important factors in disease progression.

The Lasso is a constrained version of ordinary
least squares (OLS). Let~${\bf x}_1, {\bf x}_2,\allowbreak \ldots, {\bf x}_m$ be
$n$-vectors representing the covariates, $m = 10$ and $n = 442$ in the
diabetes study, and let  ${\bf y}$ be the vector of responses for the $n$
cases. By location and scale transformations we can always assume that
the covariates have been standardized to have mean 0 and unit length,
and that the response has mean 0,
\begin{equation}\label{eqn1.1}
\quad \sum_{i=1}^n y_i = 0, \qquad \sum_{i=1}^n x_{ij} = 0, \qquad
\sum_{i=1}^n x_{ij}^2 = 1 \qquad \mbox{for }
j = 1, 2,\ldots,m. 
\end{equation}
 This is assumed to be the case in the theory which follows,
except that numerical results are expressed in the original units of
the diabetes example.

A candidate vector of regression coefficients $\widehat \bbeta =
(\widehat \beta_1, \widehat \beta_2, \ldots, \widehat \beta_m)^\prime$
gives prediction vector $\widehat \bmu$,
\begin{equation}\label{eqn1.2}
\widehat \bmu = \sum_{j=1}^m {\bf x}_j \widehat \beta_j
= X \widehat \bbeta \qquad  [X_{n \times m} = ({\bf x}_1, {\bf
  x}_2, \ldots, {\bf x}_m)]
\end{equation}
with total squared error
\begin{equation}\label{eqn1.3}
S(\widehat \bbeta) = \Vert {\bf y} - \widehat \bmu \Vert^2 =
\sum_{i=1}^n (y_i - \widehat \mu_i)^2.
\end{equation}
Let $T(\widehat \bbeta)$ be the absolute norm of $\widehat \bbeta$, 
\begin{equation}\label{eqn1.4}
T(\widehat \bbeta) = \sum_{j=1}^m \vert \widehat \beta_j \vert.
\end{equation}
The Lasso chooses $\widehat \bbeta$ by minimizing $S(\widehat
\bbeta)$ subject to a bound $t$ on $T(\widehat \bbeta)$,
\begin{equation} \label{eqn1.5}
\mbox{{\it Lasso}:}\quad \mbox{minimize} \quad S(\widehat \bbeta) \qquad
\mbox{subject to} \quad T(\widehat \bbeta) \leq t.
\end{equation}
Quadratic programming techniques can be used to solve (\ref{eqn1.5})
though we will present an easier method here, closely related to the
``homotopy method'' of \citet{osborne00}.

\begin{figure}[t]

\includegraphics{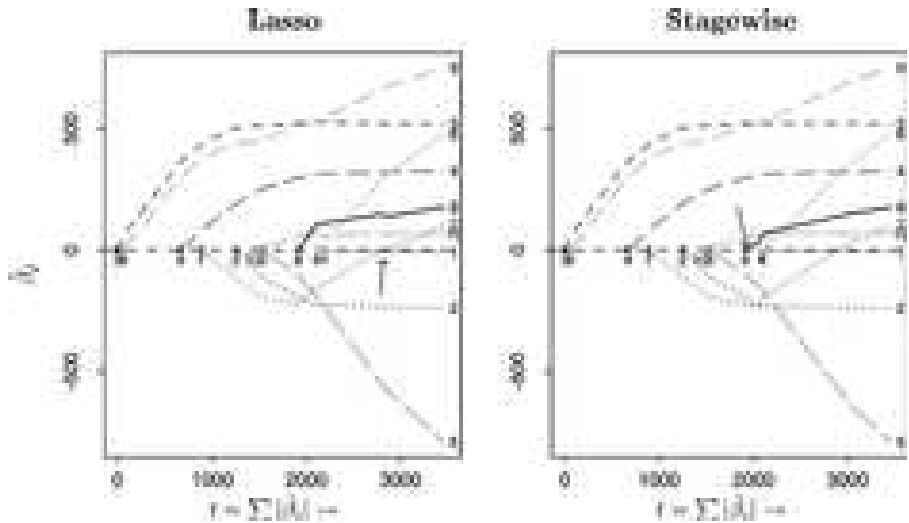}

 \caption{Estimates of
    regression coefficients $\widehat \beta_j$, $j = 1, 2, \ldots, 10$,
    for the diabetes study. ({\rm Left panel}) Lasso estimates, as a
    function of $t = \sum_j \vert \widehat \beta_j \vert$. The
    covariates enter the regression equation sequentially as $t$
    increases, in order $j = 3, 9, 4, 7, \ldots,1$. ({\rm Right panel})
    The same plot for Forward Stagewise Linear Regression. The two
    plots are nearly identical, but differ slightly for large $t$ as
    shown in the track of covariate $8$.}
\label{fig1}
\end{figure}

The left panel of Figure~\ref{fig1} shows all Lasso solutions $\widehat
\bbeta(t)$ for the diabetes study, as $t$ increases from 0, where
$\widehat \bbeta = 0$, to $t = 3460.00$, where $\widehat \bbeta$
equals the OLS regression vector, the constraint in (\ref{eqn1.5}) no longer
binding. We see that the Lasso tends to shrink the OLS
coefficients toward 0, more so for small values of $t$. Shrinkage
often improves prediction accuracy, trading off decreased variance for
increased bias as discussed in \citet{ElemStatLearn}.

The Lasso also has a parsimony property: for any given 
constraint
value $t$, only a subset of the covariates have nonzero values of
$\widehat \beta_j$. At $t = 1000$, for example, only variables 3, 9, 4
and 7 enter the Lasso regression model (\ref{eqn1.2}). If this model provides
adequate predictions, a crucial question considered in Section~\ref{sec:sec4}, the
statisticians could report these four variables as the important
ones. 

Forward Stagewise Linear Regression, henceforth called {\it
Stagewise}, is an iterative technique that begins with $\widehat
\bmu = 0$ and builds up the regression function in successive small
steps. If $\widehat \bmu$ is the current Stagewise estimate, let ${\bf
c}(\widehat \bmu)$ be the vector of {\it current correlations}
\begin{equation}\label{eqn1.6}
\widehat {\bf c} = {\bf c} (\widehat \bmu) = X^\prime ({\bf y} -
\widehat \bmu),
\end{equation}
so that $\widehat c_j$ is proportional to the correlation between
covariate $x_j$ and the current residual vector. The next step of the
Stagewise algorithm is taken in the direction of the greatest current
correlation, 
\begin{equation}\label{eqn1.7}
\widehat j = \arg\max  \vert \widehat c_j \vert \quad \mbox{and}
\quad \widehat \bmu \rightarrow \widehat \bmu + \varepsilon \cdot
\sign  (\widehat c_{\hat j}) \cdot {\bf x}_{\hat j},
\end{equation}
with $\varepsilon$ some small constant.
``Small'' is important here: the ``big'' choice $\varepsilon = \vert
\widehat c_{\hat j} \vert$ leads to the classic Forward Selection
technique, which can be overly greedy,
impulsively eliminating covariates which are correlated with $x_{\hat
  j}$. The Stagewise procedure is related to boosting and also to
Friedman's MART algorithm [\citet{Fr99}]; 
see Section~\ref{sec:sec8}, as
well as \citeauthor{ElemStatLearn} [(\citeyear{ElemStatLearn}), 
Chapter~10 and Algorithm~10.4].

The right panel of Figure~\ref{fig1} shows the coefficient plot for Stagewise
applied to the diabetes data. The estimates were built up in 6000
Stagewise steps [making $\varepsilon$~in (\ref{eqn1.7}) small enough to conceal
the ``Etch-a-Sketch'' staircase seen in Figure~\ref{fig2}, Section~2]. 
The striking fact
is the similarity between the Lasso and Stagewise estimates. Although
their definitions look completely different, the results are nearly,
{\it but not exactly}, identical.

The main point of this paper is that both Lasso and Stagewise are
variants of a basic procedure called Least Angle Regression,
abbreviated LARS (the ``S'' suggesting ``Lasso'' and ``Stagewise'').
Section~\ref{sec:sec2} describes the LARS algorithm while Section~\ref{sec:sec3} discusses
modifications that turn LARS into Lasso or Stagewise, reducing the
computational burden by at least an order of magnitude for either one.
Sections 5 and 6 verify the connections stated in Section~\ref{sec:sec3}.

Least Angle Regression is interesting in its own right, its simple
structure lending itself to inferential analysis. Section~\ref{sec:sec4} 
analyzes
the ``degrees of freedom'' of a LARS regression estimate. This leads
to a $C_p$ type statistic that suggests which estimate we should
prefer among a collection of possibilities like those in Figure~\ref{fig1}. A
particularly simple $C_p$ approximation, requiring no additional
computation beyond that for the $\widehat \bbeta$ vectors, is available
for LARS.

Section~\ref{sec:sec7} briefly discusses computational questions. An efficient $S$
program for all three methods, LARS, Lasso and Stagewise, is
available. Section~\ref{sec:sec8} elaborates on the connections with boosting.

\section{The LARS algorithm.}\label{sec:sec2}  
Least Angle Regression is a
stylized version of~the Stagewise procedure that uses a simple
mathematical formula to accelerate the computations. Only $m$ steps
are required for the full set of solutions, where $m$ is the number of
covariates: $m = 10$ in the diabetes example compared to the 6000
steps used in the right panel of Figure~\ref{fig1}. This section
describes the LARS algorithm. Modifications of LARS that produce Lasso
and Stagewise solutions are discussed in Section~\ref{sec:sec3}, and
verified in Sections~\ref{sec:sec5} and \ref{sec:sec6}.
Section~\ref{sec:sec4} uses the simple structure of LARS to help
analyze its estimation properties.

The LARS procedure works roughly as follows. As with classic Forward
Selection, we start with all coefficients equal to zero, and find the
predictor most correlated with the response, say $x_{j_1}$. We take
the largest step possible in the direction of this predictor until
some other predictor, say $x_{j_2}$, has as much correlation with the
current residual. At this point LARS parts company with Forward
Selection. Instead of continuing along $x_{j_1}$, LARS proceeds in a
direction equiangular between the two predictors until a third
variable $x_{j_3}$ earns its way into the ``most correlated'' set.
LARS then proceeds equiangularly between $x_{j_1}, x_{j_2}$ 
and~$x_{j_3}$, that is, along the ``least angle direction,'' until a fourth
variable enters, and so~on.

The remainder of this section describes the algebra necessary to
execute the equiangular strategy. As usual the algebraic details look
more complicated than the simple underlying geometry, but they lead to
the highly efficient computational algorithm described 
in Section~\ref{sec:sec7}.

\begin{figure}[b]

\includegraphics[scale=1.1]{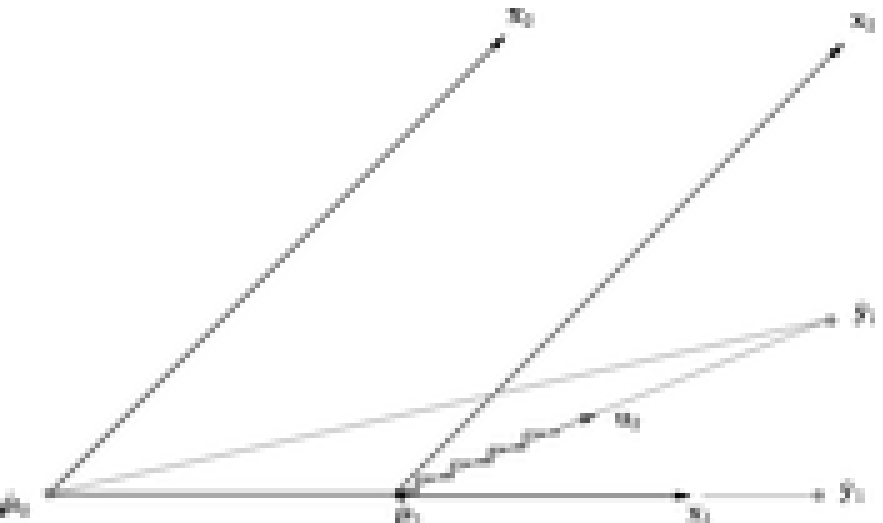}

\caption{The LARS algorithm in
the case of $m = 2$ covariates; 
$\bar {\bf y}_2$ is the projection
of ${\bf y}$ into ${\mathcal L}({\bf x}_1, {\bf x}_2)$. 
Beginning at
    $\widehat \bmu_0 = {\bf 0}$, the residual vector $\bar {\bf y}_2 -
    \widehat \bmu_0$ has greater correlation with ${\bf x}_1$ than
    ${\bf x}_2$; the next LARS estimate is $\widehat \bmu_1 = \widehat
    \bmu_0 + \widehat \gamma_1 {\bf x}_1$, where $\widehat \gamma_1$
    is chosen such that $\bar {\bf y}_2 - \widehat \bmu_1$ bisects the
    angle between ${\bf x}_1$ and ${\bf x}_2$; then $\widehat \bmu_2 =
    \widehat \bmu_1 + \widehat \gamma_2 {\bf u}_2$, where ${\bf u}_2$ is
    the unit bisector; $\widehat \bmu_2 = \bar {\bf y}_2$ in the case
    $m = 2$, but not for the case $m>2$; see Figure~\textup{\ref{fig4}}.
   The staircase indicates a typical Stagewise path. Here
    LARS gives the Stagewise track as $\varepsilon \rightarrow 0$, but a
    modification is necessary to guarantee agreement in higher
    dimensions; see Section~\textup{\ref{sec:sec3.2}}.}
\label{fig2}
\end{figure}

LARS builds up estimates $\widehat \bmu = X \widehat
\bbeta$, (\ref{eqn1.2}), in successive steps, each step adding one covariate to
the model, so that after $k$ steps just $k$ of the $\widehat
\beta_j$'s are nonzero. 
Figure~\ref{fig2} illustrates the algorithm in the
situation with $m = 2$ covariates, $X = ({\bf x}_1, {\bf x}_2)$. In
this case the current correlations (\ref{eqn1.6}) depend only on the 
projection~$\bar {\bf y}_2$ of ${\bf y}$ into the linear space ${\mathcal L}(X)$
spanned by ${\bf x}_1$ and ${\bf x}_2$,
\begin{equation} \label{eqn2.1}
{\bf c}(\widehat \bmu) = X^\prime ({\bf y} - \widehat \bmu) = X^\prime
(\bar {\bf y}_2 - \widehat \bmu).
\end{equation}

The algorithm begins at $\widehat \bmu_0 = {\bf 0}$ 
[remembering that
the response has had its mean subtracted off, as in (\ref{eqn1.1})].
Figure~\ref{fig2}
has $\bar {\bf y}_2 - \widehat \bmu_0$ making a smaller angle with
${\bf x}_1$ than ${\bf x}_2$, that is, $c_1(\widehat \bmu_0) >
c_2(\widehat \bmu_0)$. LARS then augments $\widehat \bmu_0$ in the
direction of ${\bf x}_1$, to
\begin{equation}\label{eqn2.2}
\widehat \bmu_1 = \widehat \bmu_0 + \widehat \gamma_1 {\bf x}_1.
\end{equation}
Stagewise would choose $\widehat \gamma_1$ equal to some small
value $\varepsilon$, and then repeat the process many times. 
Classic Forward Selection would take $\widehat \gamma_1$ large enough
to make $\widehat \bmu_1$ equal $\bar {\bf y}_1$, the projection of
${\bf y}$ into ${\mathcal L}({\bf x}_1)$. LARS uses an intermediate value
of $\widehat \gamma_1$, the value that makes $\bar{\bf y}_2 - \widehat
\bmu$, {\it equally} correlated with ${\bf x}_1$ and ${\bf x}_2$; that
is, $\bar {\bf y}_2 - \widehat \bmu_1$ bisects the angle between ${\bf
x}_1$ and ${\bf x}_2$, so $c_1 (\widehat \bmu_1) = c_2(\widehat
\bmu_1)$.

Let ${\bf u}_2$ be the unit vector lying along the bisector. The next
LARS estimate is 
\begin{equation}\label{eqn2.3}
\widehat \bmu_2 = \widehat \bmu_1 + \widehat \gamma_2 {\bf u}_2, 
\end{equation}
with $\widehat \gamma_2$ chosen to make $\widehat \bmu_2 = \bar
{\bf y}_2$ in the case $m = 2$. With $m > 2$ covariates, 
$\widehat\gamma_2$~would be smaller, leading to another change of direction, as
illustrated in Figure~\ref{fig4}. 
The ``staircase'' in Figure~\ref{fig2} indicates a
typical Stagewise path. LARS is motivated by the fact that it is easy
to calculate the step sizes $\widehat \gamma_1, \widehat \gamma_2,
\ldots$ theoretically, short-circuiting the small Stagewise steps.

Subsequent LARS steps, beyond two covariates, are taken along {\it
  equiangular vectors}, generalizing the bisector ${\bf u}_2$ in
Figure~\ref{fig2}. 
{\it We assume that the covariate vectors ${\bf x}_1, {\bf
    x}_2, \ldots, {\bf x}_m$ are linearly independent}. For ${\mathcal A}$
a subset of the indices $\{1, 2, \ldots, m\}$, define the matrix
\begin{equation} \label{eqn2.4}
X_{\mathcal A} = (\cdots \,s_j {\bf x}_j\, \cdots)_{j \in {\mathcal A}},
\end{equation}
where the signs $s_j$ equal $\pm 1$. Let
\begin{equation}\label{eqn2.5}
{\mathcal G}_{\mathcal A} = X_{\mathcal A}^\prime X_{\mathcal A} 
\quad \mbox{and}
\quad A_{\mathcal A} = (\bolds{1}_{\mathcal A}^\prime 
{\mathcal G}_{\mathcal A}^{-1}
\bolds{1}_{\mathcal A})^{-1/2},
\end{equation}
$\bolds{1}_{\mathcal A}$ being a vector of 1's of length equaling $|\mathcal A|$, the size of
${\mathcal A}$. The
\begin{equation}\label{eqn2.6}
\mbox{\it equiangular vector}\quad  
{\bf u}_{\mathcal A} = X_{\mathcal A}w_{\mathcal A}
\qquad \mbox{where } w_{\mathcal A} = A_{\mathcal A} G_{\mathcal A}^{-1}
\bolds{1}_{\mathcal A},
\end{equation}
is the unit vector making equal angles, less than $90^{\circ}$,
with the columns of $X_{\mathcal A}$,
\begin{equation} \label{eqn2.7}
X_{\mathcal A}^\prime {\bf u}_{\mathcal A} = A_{\mathcal A} 
\bolds{1}_{\mathcal A} \quad
\mbox{and} \quad \Vert {\bf u}_{\mathcal A} \Vert^2 = 1.
\end{equation}

We can now fully describe the LARS algorithm. As with the Stagewise
procedure we begin at $\widehat \bmu_0 = {\bf 0}$ and build up
$\widehat \bmu$ by steps, larger steps in the LARS case. Suppose that
$\widehat \bmu_{\mathcal A}$ is the current LARS estimate and that 
\begin{equation}\label{eqn2.8}
\widehat {\bf c} = X^\prime ({\bf y} - \widehat \bmu_{\mathcal A})
\end{equation}
is the vector of current correlations (\ref{eqn1.6}). The {\it active set}
${\mathcal A}$ is the set of indices corresponding to covariates with the
greatest absolute current correlations,
\begin{equation}\label{eqn2.9}
\widehat C = \max_{j} \{\vert \widehat c_j \vert \} \quad \mbox{and}
\quad {\mathcal A} = \{j \dvtx \vert \widehat c_j \vert = \widehat C\}.
\end{equation}

Letting
\begin{equation}\label{eqn2.10}
s_j = \sign \{\widehat c_j\} \qquad \mbox{for }  j \in {\mathcal
  A},
\end{equation}
we compute $X_{\mathcal A}, A_{\mathcal A}$ and ${\bf u}_{\mathcal A}$ as in
(\ref{eqn2.4})--(\ref{eqn2.6}), and also the inner product vector
\begin{equation} \label{eqn2.11}
{\bf a} \equiv X^\prime {\bf u}_{\mathcal A}.
\end{equation}
Then the next step of the LARS algorithm updates $\widehat
\bmu_{\mathcal A}$, say to
\begin{equation}\label{eqn2.12}
\widehat \bmu_{{\mathcal A}_+} = \widehat \bmu_{\mathcal A} + \widehat \gamma
{\bf u}_{\mathcal A},
\end{equation}
where
\begin{equation}\label{eqn2.13}
\widehat \gamma = \min_{j \in {\mathcal A}^c}{}^{\!\! +}
\biggl\{\frac{\widehat C - \widehat c_j}{A_{\mathcal A} - a_j},
\frac{\widehat C + \widehat c_j}{A_{\mathcal A} + a_j} \biggr\};
\end{equation}
``$\min^+$'' 
indicates that the minimum is taken over only
positive components within each choice of $j$ in (\ref{eqn2.13}).

\begin{figure}[b]

\includegraphics[scale=0.9]{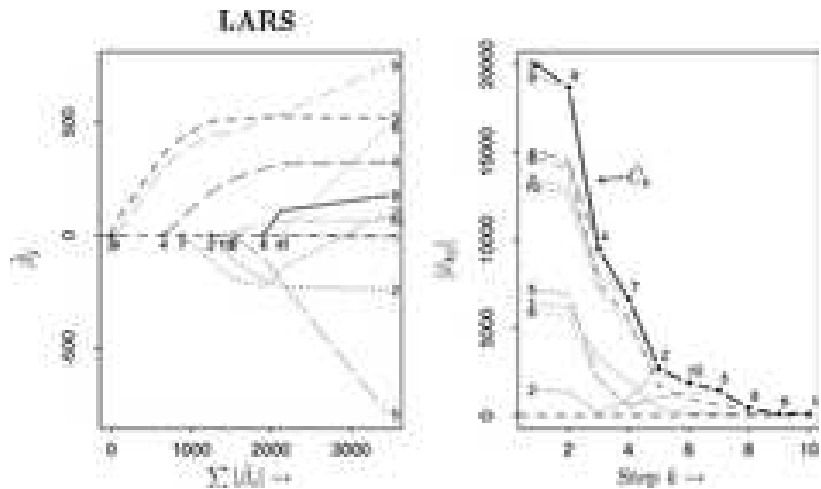}

\caption{LARS analysis of the
    diabetes study: ({\rm left}) estimates of regression 
coefficients~$\widehat \beta_j$,  $j = 1, 2, \ldots, 10$; plotted 
versus $\sum
    \vert \widehat \beta_j \vert$; plot is slightly different than
    either Lasso or Stagewise, Figure~\textup{\ref{fig1}}; ({\rm right}) 
absolute current
    correlations as function of LARS step; variables enter active 
set~$(\ref{eqn2.9})$ in order $3, 9, 4, 7, \ldots, 1$; 
heavy curve shows maximum
current correlation $\widehat C_k$ declining with $k$.}
\label{fig3}
\end{figure}

Formulas (\ref{eqn2.12}) and (\ref{eqn2.13}) have the following interpretation: define
\begin{equation}\label{eqn2.14}
\bmu (\gamma) = \widehat \bmu_{\mathcal A} + \gamma {\bf u}_{\mathcal A},
\end{equation}
for $\gamma > 0$, so that the current correlation
\begin{equation}\label{eqn2.15}
c_j(\gamma) = {\bf x}_j^\prime \bigl({\bf y} - \bmu(\gamma)\bigr) = \widehat c_j
- \gamma a_j.
\end{equation}
For $j \in {\mathcal A}$, (\ref{eqn2.7})--(\ref{eqn2.9}) yield
\begin{equation}\label{eqn2.16}
\vert c_j (\gamma) \vert = \widehat C - \gamma A_{\mathcal A},
\end{equation}
showing that all of the maximal absolute current correlations
decline equally. For $j \in {\mathcal A}^c$, equating (\ref{eqn2.15}) with (\ref{eqn2.16})
shows that $c_j(\gamma)$ equals the maximal value at $\gamma =
(\widehat C - \widehat c_j) / (A_{\mathcal A} - a_j)$. Likewise
$-c_j(\gamma)$, the current correlation for the reversed covariate
$-{\bf x}_j$, achieves maximality at $(\widehat C + \widehat c_j) /
(A_{\mathcal A} + a_j)$. Therefore {\it $\widehat \gamma$ in {\rm (\ref{eqn2.13})}
  is the smallest positive value of $\gamma$ such that some new index
  $\widehat j$ joins the active set}; $\widehat j$ is the minimizing
index in (\ref{eqn2.13}), and the new active set ${\mathcal A}_+$ is ${\mathcal A} \cup
\{ \widehat j\}$; the new maximum absolute correlation is $\widehat
C_+ = \widehat C - \widehat \gamma A_{\mathcal A}$.

Figure~\ref{fig3} concerns the LARS analysis of the diabetes data. The complete
algorithm required\ only $m = 10$ steps of procedure (\ref{eqn2.8})--(\ref{eqn2.13}), with
the variables joining the active set ${\mathcal A}$ in the same order as
for the Lasso: $3, 9, 4, 7, \ldots, 1$. Tracks of the regression
coefficients $\widehat \beta_j$ are nearly but not exactly the same as
either the Lasso or Stagewise tracks of Figure~\ref{fig1}.

The right panel shows the absolute current correlations 
\begin{equation}\label{eqn2.17}
\vert \widehat c_{kj} \vert = \vert {\bf x}_j^\prime ({\bf y} -
\widehat \bmu_{k-1})\vert
\end{equation}
for variables $j = 1, 2, \ldots, 10$, as a function of the LARS
step $k$. The maximum correlation
\begin{equation}\label{eqn2.18}
\widehat C_k = \max\{\vert \widehat c_{kj} \vert \} = \widehat C_{k-1}
- \widehat \gamma_{k-1} A_{k-1}
\end{equation}
declines with $k$, as it must. At each step a new variable $j$
joins the active set, henceforth having $\vert \widehat c_{kj} \vert =
\widehat C_k$. The sign $s_j$ of each ${\bf x}_j$ in (\ref{eqn2.4}) stays
constant as the active set increases.

\begin{figure}[b]

\includegraphics[scale=1.1]{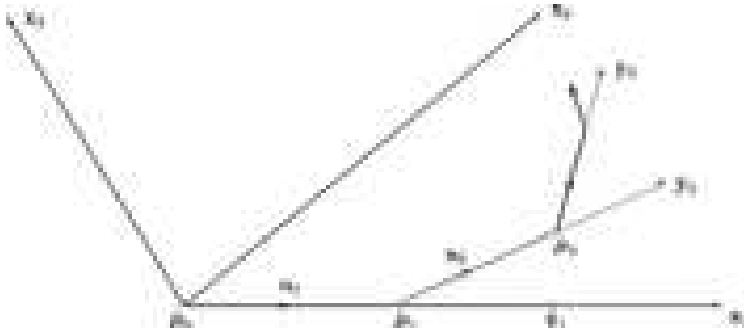}

\caption{At each stage the LARS
  estimate $\widehat \bmu_k$ approaches, but does not reach, the
  corresponding OLS estimate $\bar {\bf y}_k$.}
\label{fig4}
\end{figure}

Section~\ref{sec:sec4} makes use of the relationship between Least Angle Regression
and Ordinary Least Squares illustrated in Figure~\ref{fig4}. Suppose LARS has
just completed step $k-1$, giving $\widehat \bmu_{k-1}$, and is
embarking upon step $k$. The active set ${\mathcal A}_k$, (\ref{eqn2.9}), will have
$k$ members, giving $X_k, {\mathcal G}_k, A_k$ and ${\bf u}_k$ as in
(\ref{eqn2.4})--(\ref{eqn2.6}) (here replacing subscript ${\mathcal A}$ with ``$k$''). Let
$\bar {\bf y}_k$ indicate the projection of ${\bf y}$ into ${\mathcal
  L}(X_k)$, which, since $\widehat \bmu_{k-1} \in {\mathcal L}(X_{k-1})$,
is
\begin{equation}\label{eqn2.19}
\bar {\bf y}_k = \widehat \bmu_{k-1} + X_k {\mathcal G}_k^{-1} X_k^\prime
({\bf y} - \widehat \bmu_{k-1}) = \widehat \bmu_{k-1} + \frac{\widehat
  C_k}{A_k} \ {\bf u}_k,
\end{equation}
the last equality following from (\ref{eqn2.6}) and the fact that the
signed current correlations in ${\mathcal A}_k$ all equal $\widehat C_k$,
\begin{equation} \label{eqn2.20}
X_k^\prime ({\bf y} - \widehat \bmu_{k-1}) = \widehat C_k \bolds{1}_{\mathcal A}.
\end{equation}

Since ${\bf u}_k$ is a unit vector, (\ref{eqn2.19}) says that $\bar {\bf y}_k -
\widehat \bmu_{k-1}$ has length
\begin{equation}\label{eqn2.21}
\bar \gamma_k \equiv \frac{\widehat C_k}{A_k}.
\end{equation}

Comparison with\ (\ref{eqn2.12}) shows that the LARS estimate $\widehat
\bmu_k$ lies on the line from~$\widehat \bmu_{k-1}$ to $\bar {\bf
  y}_k$,
\begin{equation}\label{eqn2.22}
\widehat \bmu_k - \widehat \bmu_{k-1} = \frac{\widehat \gamma_k}{\bar
  \gamma_k} \ (\bar {\bf y}_k - \widehat \bmu_{k-1}).
\end{equation}
It is easy to see that $\widehat \gamma_k$, (\ref{eqn2.12}), is always
less than $\bar \gamma_k$, so that $\widehat \bmu_k$ lies closer than
$\bar {\bf y}_k$ to $\widehat \bmu_{k-1}$. Figure~\ref{fig4} shows the
successive LARS estimates $\widehat \bmu_k$ always approaching but
never reaching the OLS estimates $\bar {\bf y}_k$. 

The exception is at the last stage: since ${\mathcal A}_m$ contains all
covariates, (\ref{eqn2.13}) is not defined. By convention the algorithm takes
$\widehat \gamma_m = \bar \gamma_m = \widehat C_m / A_m$, making
$\widehat \bmu_m = \bar {\bf y}_m$ and $\widehat \bbeta_m$ equal the
OLS estimate for the full set of $m$ covariates.

The LARS algorithm is computationally thrifty. Organizing the
calculations correctly, the computational cost for the
entire $m$ steps is of the same order as that required  for the usual Least
Squares solution for the full set of $m$ covariates. Section~\ref{sec:sec7}
describes an efficient LARS program available from the authors. With
the modifications described in the next section, this program also
provides economical Lasso and Stagewise solutions.

\section{Modified versions of Least Angle Regression.}\label{sec:sec3}

Figures 1 and 3 show Lasso, Stagewise and LARS yielding remarkably
similar estimates for the diabetes data. The similarity is no
coincidence. This section describes simple modifications of the LARS
algorithm that produce Lasso or Stagewise estimates. Besides improved
computational efficiency, these relationships elucidate the methods'
rationale: all three algorithms can be viewed as moderately greedy
forward stepwise procedures whose forward progress is determined by
compromise among the currently most correlated covariates. LARS moves
along the most obvious compromise direction, the equiangular vector
(\ref{eqn2.6}), while Lasso and Stagewise put some restrictions on the
equiangular strategy.

\subsection{The LARS--Lasso relationship.}\label{sec:sec3.1}
The full set of Lasso solutions, as shown for the diabetes study in
Figure~\ref{fig1}, can be generated by a minor modification of the
LARS algorithm (\ref{eqn2.8})--(\ref{eqn2.13}). Our main result is
described here and verified in Section~\ref{sec:sec5}. It closely
parallels the homotopy method in the papers by Osborne, Presnell and
Turlach (\citeyear{osborne00}, b),
though the LARS approach is somewhat more direct.

Let $\widehat \bbeta$ be a Lasso solution (\ref{eqn1.5}), with $\widehat \bmu =
X \widehat \bbeta$. Then it is easy to show that the sign of any
nonzero coordinate $\widehat \beta_j$ must agree with the sign $s_j$
of the current correlation $\widehat c_j = {\bf x}_j^\prime ({\bf y} -
\widehat \bmu)$,
\begin{equation}\label{eqn3.1}
\sign (\widehat \beta_j) = \sign (\,\widehat c_j) = s_j;
\end{equation}
see Lemma~\ref{lem:lem5} of Section~\ref{sec:sec5}. The LARS algorithm does not enforce
restriction (\ref{eqn3.1}), but it can easily be modified to do so.

Suppose we have just completed a LARS step, giving a new active set
${\mathcal A}$ as in (\ref{eqn2.9}), and that the corresponding LARS estimate
$\widehat \bmu_{\mathcal A}$ corresponds to a Lasso solution $\widehat
\bmu = X \widehat \bbeta$. Let
\begin{equation}\label{eqn3.2}
w_{\mathcal A} = A_{\mathcal A} {\mathcal G}_{\mathcal A}^{-1} \bolds{1}_{\mathcal A},
\end{equation}
a vector of length the size of ${\mathcal A}$, and (somewhat abusing
subscript notation) define~$\widehat {\bf d}$ to be the $m$-vector
equaling $s_j w_{{\mathcal A} j}$ for $j \in {\mathcal A}$ and zero
elsewhere. Moving in the positive $\gamma$ direction along the LARS
line (\ref{eqn2.14}), we see that
\begin{equation}\label{eqn3.3}
\bmu (\gamma) = X \bbeta(\gamma), \qquad \mbox{where }  \beta_j
(\gamma) = \widehat \beta_j + \gamma \widehat d_j
\end{equation}
for $j \in {\mathcal A}$. Therefore $\beta_j(\gamma)$ will change
sign at 
\begin{equation}\label{eqn3.4}
\gamma_j = - \widehat \beta_j / \widehat d_j,
\end{equation}
the first such change occurring at
\begin{equation}\label{eqn3.5}
\widetilde \gamma = \min_{\gamma_j > 0} \{ \gamma_j \},
\end{equation}
say for covariate $x_{\tilde j}$; $\widetilde \gamma$ equals
infinity by definition if there is no $\gamma_j > 0$.

If $\widetilde \gamma$ is less than $\widehat \gamma$, (\ref{eqn2.13}), then
$\beta_j(\gamma)$ cannot be a Lasso solution for $\gamma > \widetilde
\gamma$ since the sign restriction (\ref{eqn3.1}) must be violated:
$\beta_{\tilde j}(\gamma)$ has changed sign while $c_{\tilde
  j}(\gamma)$ has not. [The continuous function $c_{\tilde j}(\gamma)$
cannot change sign within a single LARS step since $\vert c_{\tilde
  j}(\gamma) \vert = \widehat C - \gamma A_{\mathcal A} > 0$, 
(\ref{eqn2.16}).]

\begin{las*} 
If $\widetilde \gamma < \widehat
\gamma$, stop the ongoing LARS step at $\gamma = \widetilde \gamma$
and remove $\widetilde j$ from the calculation of the next equiangular
direction. That is,
\begin{equation}\label{eqn3.6}
\widehat \bmu_{{\mathcal A}_+} = 
\widehat \bmu_{\mathcal A} + \widetilde
\gamma {\bf u}_{\mathcal A} \quad \mbox{and} \quad 
{\mathcal A}_+ = {\mathcal A} -
\{\widetilde j\}
\end{equation}
rather than (\ref{eqn2.12}).
\end{las*}

\begin{thm}
\label{thm:thm1}
 Under the Lasso modification, and assuming
the \textup{``}one at a time\textup{''} condition discussed below, the LARS algorithm
yields all Lasso solutions.
\end{thm}

The active sets ${\mathcal A}$ grow monotonically larger as the original
LARS algorithm progresses, but the Lasso modification allows ${\mathcal
  A}$ to decrease. ``One at a time'' means that the increases and
decreases never involve more than a single index~$j$. This is the
usual case for quantitative data and can always be realized by adding
a little jitter to the $y$ values. Section~\ref{sec:sec5} discusses tied
situations.

The Lasso diagram in Figure~\ref{fig1} was actually calculated using the
modified LARS algorithm. Modification (\ref{eqn3.6}) came into play only once,
at the arrowed point in the left panel. There ${\mathcal A}$ contained all
10 indices while ${\mathcal A}_+ = {\mathcal A} - \{7\}$. Variable~7 was
restored to the active set one LARS step later, the next and last step
then taking $\widehat \bbeta$ all the way to the full OLS solution. The
brief absence of variable 7 had an effect on the tracks of the others,
noticeably $\widehat \beta_8$. The price of using Lasso instead of
unmodified LARS comes in the form of added steps, 12 instead of 10 in
this example. For the more complicated ``quadratic model'' of Section
4, the comparison was 103 Lasso steps versus 64 for LARS.

\subsection{The LARS--Stagewise relationship.}\label{sec:sec3.2}

The staircase in Figure~\ref{fig2} indicates how the Stagewise algorithm might
proceed forward from $\widehat \bmu_1$, a point of equal current
correlations $\widehat c_1 = \widehat c_2$, 
(\ref{eqn2.8}). The first small
step has (randomly) selected index $j=1$, taking us to $\widehat
\bmu_1 + \varepsilon {\bf x}_1$. Now variable 2 is more correlated,
\begin{equation}\label{eqn3.7}
{\bf x}_2^\prime ({\bf y} - \widehat \bmu_1 - \varepsilon {\bf x}_1) >
{\bf x}_1^\prime ({\bf y} - \widehat \bmu_1 - \varepsilon {\bf x}_1),
\end{equation}
forcing $j = 2$ to be the next Stagewise choice and so on.

{\it We will consider an idealized Stagewise procedure in which the
  step size $\varepsilon$~goes to zero.} This collapses the staircase
along the direction of the bisector ${\bf u}_2$ in Figure~\ref{fig2}, making
the Stagewise and LARS estimates agree. They always agree for $m = 2$
covariates, but another modification is necessary for LARS to produce
Stagewise estimates in general. Section~\ref{sec:sec6} verifies the main result
described next.

Suppose that the Stagewise procedure has taken $N$ steps of
infinitesimal size $\varepsilon$~from some previous estimate $\widehat
\bmu$, with
\begin{equation}\label{eqn3.8}
N_j \equiv \# \{ \mbox{steps with selected index } j\}, \qquad j = 1,
2, \ldots, m.
\end{equation}
It is easy to show, as in Lemma~\ref{lem:lem8} of Section~\ref{sec:sec6}, that $N_j = 0$
for $j$ not in the active set ${\mathcal A}$ defined by the current
correlations ${\bf x}_j^\prime ({\bf y} - \widehat \bmu)$, (\ref{eqn2.9}).
Letting
\begin{equation}\label{eqn3.9}
P \equiv (N_1, N_2, \ldots, N_m) / N,
\end{equation}
with $P_{\mathcal A}$ indicating the coordinates of $P$ for $j \in
{\mathcal A}$, the new estimate is
\begin{equation}\label{eqn3.10}
\bmu = \widehat \bmu + N \varepsilon X_{\mathcal A} P_{\mathcal A}\qquad 
[(\ref{eqn2.4})].
\end{equation}
(Notice that the Stagewise steps are taken along the
directions $s_j {\bf x}_j$.)

The LARS algorithm (\ref{eqn2.14}) progresses along
\begin{equation}\label{eqn3.11}
\bmu_{\mathcal A} + \gamma X_{\mathcal A} w_{\mathcal A}, \qquad 
\mbox{where } 
w_{\mathcal A} = A_{\mathcal A} {\mathcal G}_{\mathcal A}^{-1} 
\bolds{1}_{\mathcal A}\qquad [(\ref{eqn2.6})\mbox{--}(\ref{eqn3.2})].
\end{equation}
Comparing (\ref{eqn3.10}) with (\ref{eqn3.11}) shows that LARS cannot
agree with Stagewise if $w_{\mathcal A}$ has negative components, since
$P_{\mathcal A}$ is nonnegative. To put it another way, the direction of
Stagewise progress $X_{\mathcal A} P_{\mathcal A}$ must lie in the convex cone
generated by the columns of $X_{\mathcal A}$,
\begin{equation}\label{eqn3.12}
{\mathcal C}_{\mathcal A} = \Biggl\{{\bf v} = \sum_{j \in {\mathcal A}} s_j {\bf
 x}_j P_j, \ P_j \geq 0 \Biggr\}.
\end{equation}

If ${\bf u}_{\mathcal A} \in {\mathcal C}_{\mathcal A}$ then there is no
contradiction between (\ref{eqn3.12}) and (\ref{eqn3.13}). If not it seems natural to
replace ${\bf u}_{\mathcal A}$ with its projection into ${\mathcal C}_{\mathcal
  A}$, that is, the nearest point in the convex cone.

\begin{stag*}
Proceed as in (\ref{eqn2.8})--(\ref{eqn2.13}),
except with ${\bf u}_{\mathcal A}$ replaced by ${\bf u}_{\hat {\mathcal B}}$,
the unit vector lying along the projection of ${\bf u}_{\mathcal A}$ into
${\mathcal C}_{\mathcal A}$. (See Figure~\ref{fig8} 
in Section~\ref{sec:sec6}.)
\end{stag*}

\begin{thm}
\label{thm:thm2}  
 Under the Stagewise modification, the LARS
algorithm yields all Stagewise solutions.
\end{thm}

The vector ${\bf u}_{\hat {\mathcal B}}$ in the Stagewise modification is
the equiangular vector (\ref{eqn2.6}) for the subset $\widehat {\mathcal B}
\subseteq {\mathcal A}$ corresponding to the face of ${\mathcal C}_{\mathcal A}$
into which the projection falls. Stagewise is a LARS type algorithm
that allows the active set to decrease by one or more indices. This
happened at the arrowed point in the right panel of Figure~\ref{fig1}: there
the set ${\mathcal A} = \{3, 9, 4, 7, 2, 10, 5, 8 \}$ was decreased to
$\widehat {\mathcal B} = {\mathcal A} - \{3, 7 \}$. It took a total of 13
modified LARS steps to reach the full OLS solution $\bar {\bbeta}_m =
(X^\prime X)^{-1} X^\prime {\bf y}$. The three methods, LARS, Lasso
and Stagewise, always reach OLS eventually, but LARS does so in only
$m$ steps while Lasso and, especially, Stagewise can 
take longer. For
the $m = 64$ quadratic model of Section~\ref{sec:sec4}, Stagewise took 255 steps.

According to Theorem~\ref{thm:thm2} the difference between successive
Stagewise--modified LARS estimates is 
\begin{equation}\label{eqn3.13}
\widehat \bmu_{{\mathcal A}_+} - \widehat \bmu_{\mathcal A} = \widehat \gamma
{\bf u}_{\hat {\mathcal B}} = \widehat \gamma X_{\hat {\mathcal B}} w_{\hat
  {\mathcal B}},
\end{equation}
as in (\ref{eqn3.13}). Since ${\bf u}_{\hat {\mathcal B}}$ exists in the
convex cone ${\mathcal C}_{\mathcal A}$, $w_{\hat {\mathcal B}}$ must have
nonnegative components.  This says that the difference of successive
coefficient estimates for coordinate $j \in {\widehat {\mathcal B}}$
satisfies
\begin{equation}\label{eqn3.14}
\sign  (\widehat \beta_{+j} - \widehat \beta_j) = s_j,
\end{equation}
where $s_j = \sign \{ {\bf x}_j^\prime ({\bf y} - \widehat
\bmu)\}$. 

We can now make a useful comparison of the three methods:
\begin{enumerate}
\item {\it Stagewise}---successive differences of $\widehat \beta_j$
  agree in sign with the current correlation $\widehat c_j = {\bf
    x}_j^\prime ({\bf y} - \widehat \bmu)$;
\item {\it Lasso}---$\widehat \beta_j$ agrees in sign with $\widehat
  c_j$;
\item {\it LARS}---no sign restrictions (but see Lemma~\ref{lem:lem1}
  of Section~\ref{sec:sec5}).
\end{enumerate}
From this point of view, Lasso is intermediate between the LARS
and Stagewise methods.

The successive difference property (\ref{eqn3.14}) makes the Stagewise
$\widehat \beta_j$ estimates move monotonically away from 0. Reversals
are possible only if $\widehat c_j$ changes sign while $\widehat
\beta_j$ is ``resting'' between two periods of change. This happened
to variable 7 in Figure~\ref{fig1} between the 8th and 10th Stagewise-modified
LARS steps.

\subsection{Simulation study.}
A small simulation study was
carried out comparing the LARS, Lasso and Stagewise algorithms. The
$X$ matrix for the simulation was based on the diabetes example of
Table~\ref{tab1}, but now using a ``Quadratic Model'' having $m = 64$
predictors, including interactions and squares of the 10 
original
covariates:
\begin{equation}
 \label{eq:insert3.1}
\mbox{\em Quadratic Model} \quad \mbox{10 main effects},
45\mbox{ interactions},  9  \mbox{ squares},  
\end{equation}
the last being the squares of each ${\bf x}_j$ except the
dichotomous variable ${\bf x}_2$. The true mean vector $\bmu$ for the
simulation was $\bmu = X \bbeta$, where $\bbeta$ was obtained by
running LARS for 10 steps on the original $(X, {\bf y})$ diabetes data
(agreeing in this case with the 10-step Lasso or Stagewise analysis).
Subtracting $\bmu$ from a centered version of the original ${\bf y}$
vector of Table~\ref{tab1} gave a vector $\bepsi = {\bf y} - \bmu$ of $n = 442$
residuals. The ``true $R^2$'' for this model, $\Vert \bmu \Vert^2 /
(\Vert \bmu \Vert^2 + \Vert \bepsi \Vert^2)$, equaled 0.416.

100 simulated response vectors ${\bf y}^\ast$ were generated from the
model
\begin{equation}
  \label{eq:insert3.2}
{\bf y}^\ast = \bmu + \bepsi^\ast,
\end{equation}
with $\bepsi^\ast = (\varepsilon_1^\ast, \varepsilon_2^\ast, \ldots,
\varepsilon_n^\ast)$ a random sample, with replacement, from the
components of $\bepsi$. The LARS algorithm with $K = 40$ steps was run
for each simulated data set $(X, {\bf y}^\ast)$, yielding a sequence
of estimates $\widehat \bmu^{(k)\ast}$, $k = 1, 2, \ldots, 40$, and
likewise using the Lasso and Stagewise algorithms.

\begin{figure}[b]

\includegraphics{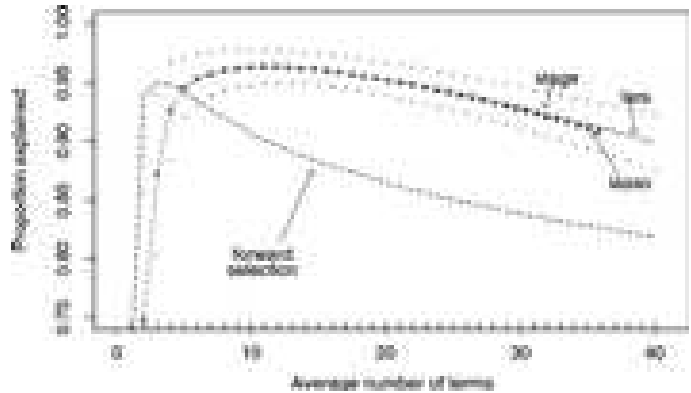}

\caption{Simulation study comparing
  LARS, Lasso and Stagewise algorithms; $100$ 
replications of model
\textup{(\ref{eq:insert3.1})}--\textup{(\ref{eq:insert3.2})}. 
Solid curve shows
  average proportion explained, \textup{(\ref{eq:insert3.3})}, for LARS
  estimates as function of number of steps $k = 1, 2, \ldots, 40$;
  Lasso and Stagewise give nearly identical results; small dots
  indicate plus or minus one standard deviation over the $100$ simulations.
  Classic Forward Selection (heavy dashed curve) rises and falls more
  abruptly.}
\label{fig:Insertfig5}
\end{figure}

Figure~\ref{fig:Insertfig5} compares the LARS, Lasso and Stagewise estimates. For a
given estimate $\widehat \bmu$ define the {\it proportion explained}
$\pe(\widehat \bmu)$ to be
\begin{equation}
  \label{eq:insert3.3}
\pe(\widehat \bmu) = 1 - \Vert \widehat \bmu - \bmu \Vert^2 / \Vert
\bmu \Vert^2,
\end{equation}
so $\pe ({\bf 0}) = 0$ and $\pe(\bmu) = 1$. The solid curve graphs
the average of $\pe(\widehat \bmu^{(k)\ast})$ over the 100 simulations,
versus step number $k$ for LARS, $k = 1, 2, \ldots, 40$. The
corresponding curves are graphed for Lasso and Stagewise, except that
the horizontal axis is now the average number of nonzero $\widehat
\beta_j^\ast$ terms composing $\widehat \bmu^{(k)\ast}$.  For example,
$\widehat \bmu^{(40)\ast}$ averaged 33.23 nonzero terms with Stagewise,
compared to 35.83 for Lasso and 40 for LARS.

Figure~\ref{fig:Insertfig5}'s most striking message is that the three algorithms
performed almost identically, and rather well. The average proportion
explained rises quickly, reaching a maximum of 0.963 at $k = 10$, and
then declines slowly as $k$ grows to 40. The light dots display the small
standard deviation of $\pe(\widehat \bmu^{(k)\ast})$ over the 100
simulations, roughly $\pm 0.02$. Stopping at any point between $k = 5$
and 25 typically gave a $\widehat \bmu^{(k)\ast}$ with true predictive
$R^2$ about $0.40$, compared to the ideal value $0.416$ for $\bmu$.

The dashed curve in Figure~\ref{fig:Insertfig5} tracks the average proportion explained
by classic Forward Selection. It rises very quickly, to a maximum of
$0.950$ after $k = 3$ steps, and then falls back more abruptly than the
LARS--Lasso--Stagewise curves. This behavior agrees with the
characterization of Forward Selection as a dangerously greedy
algorithm.

\subsection{Other LARS modifications.}\label{sec:sec3.3}
Here are a few more examples of LARS type model-building algorithms.
 
\begin{pos*}
Constraint (\ref{eqn1.5}) can be strengthened
to
\begin{equation}
\label{eqn3.15}
\quad \mbox{minimize} \quad  S(\widehat \bbeta) \qquad  
\mbox{subject to} \quad 
T(\widehat
\bbeta) \leq t \mbox{ and all }\widehat \beta_j \geq
0.
\end{equation}
\end{pos*}

This would be appropriate if the statisticians or scientists
believed that the variables $x_j$ {\it must} enter the prediction
equation in their defined directions. Situation~(\ref{eqn3.15}) is a more
difficult quadratic programming problem than (\ref{eqn1.5}), but it can be
solved by a further modification of the Lasso-modified LARS algorithm:
change $\vert \widehat c_j \vert$ to $\widehat c_j$ at both places in
(\ref{eqn2.9}), set $s_j = 1$ instead of (\ref{eqn2.10}) 
and change~(\ref{eqn2.13}) to
\begin{equation}\label{eqn3.16}
\widehat \gamma = \min_{j \in {\mathcal A}^c}{}^{\!\! +} \biggl\{
\frac{\widehat C - \widehat c_j}{A_{\mathcal A} - a_j} \biggr\}.
\end{equation}
The positive Lasso usually does {\it not} converge to the full
OLS solution $\bar \beta_m$, even for very large choices of $t$.

The changes above amount to considering the ${\bf x}_j$ as generating
half-lines rather than full one-dimensional spaces. A positive
Stagewise version can be developed in the same way, and has the
property that the $\widehat \beta_j$ tracks are always monotone.

\subsubsection*{LARS--OLS hybrid.}
After $k$ steps the LARS algorithm
has identified a set ${\mathcal A}_k$ of covariates, for example, 
${\mathcal  A}_4 = \{3, 9, 4, 7\}$ in the diabetes study. Instead of $\widehat
{\bbeta}_k$ we might prefer $\bar {\bbeta}_k$, the OLS coefficients
based on the linear model with covariates in~${\mathcal A}_k$---using
LARS to find the model but not to estimate the coefficients. Besides
looking more familiar, this will always increase the usual empirical
$R^2$ measure of fit (though not necessarily the true fitting
accuracy),
\begin{equation}
\label{eqn3.17}
R^2(\bar {\bbeta}_k) - R^2 (\widehat {\bbeta}_k) =
\frac{1-\rho^2_k}{\rho_k(2- \rho_k)} 
[R^2 (\widehat {\bbeta}_k) - R^2
(\widehat {\bbeta}_{k-1})],
\end{equation}
where $\rho_k = \widehat \gamma_k / \bar \gamma_k$ as in
(\ref{eqn2.22}). 

The increases in $R^2$ were small in the diabetes example, on the
order of 0.01 for~$k \geq 4$ compared with $R^2 \  \dot{=} \ 0.50$, which
is expected from (\ref{eqn3.17}) since we would usually continue LARS until
$R^2(\widehat {\bbeta}_k) - R^2(\widehat {\bbeta}_{k-1})$ was small.
For the same reason~$\bar {\bbeta}_k$ and $\widehat {\bbeta}_k$ are
likely to lie near each other as they did in the diabetes example.

\subsubsection*{Main effects first.} 
It is straightforward to restrict
the order in which variables are allowed to enter the LARS algorithm.
For example, having obtained ${\mathcal A}_4 = \{3, 9, 4, 7\}$ for the
diabetes study, we might {\it then} wish to check for interactions.
To do this we begin LARS again, replacing ${\bf y}$ with ${\bf y} -
\widehat \bmu_4$ and ${\bf x}$ with the $n \times 6$ matrix whose
columns represent the interactions ${\bf x}_{3:9}, {\bf x}_{3:4},
\ldots, {\bf x}_{4:7}$.

\subsubsection*{Backward Lasso.}
The Lasso--modified LARS algorithm can
be run backward, starting from the full OLS solution $\bar
{\bbeta}_m$.  Assuming that all the coordinates of $\bar {\bbeta}_m$
are nonzero, their signs must agree with the signs $s_j$ that the
current correlations had during the final LARS step. This allows us to
calculate the last equiangular direction ${\bf u}_{\mathcal A}$,
(\ref{eqn2.4})--(\ref{eqn2.6}). Moving backward 
from $\widehat \bmu_m = X \bar
{\bbeta}_m$ along the line $\bmu (\gamma) = \widehat \bmu_m - \gamma
{\bf u}_{\mathcal A}$, we eliminate from the active set the index of the
first $\widehat \beta_j$ that becomes zero. Continuing backward, we
keep track of all coefficients $\widehat \beta_j$ and current
correlations $\widehat c_j$, following essentially the same rules for
changing ${\mathcal A}$ as in Section~\ref{eqn3.1}. 
As in (\ref{eqn2.3}), (\ref{eqn3.5}) the
calculation of $\widetilde \gamma$ and $\widehat \gamma$ is easy.

The crucial property of the Lasso that makes backward navigation
possible is~(\ref{eqn3.1}), which permits calculation of the correct
equiangular direction ${\bf u}_{\mathcal A}$ at each step. In this sense
Lasso can be just as well thought of as a backward-moving algorithm.
This is not the case for LARS or Stagewise, both of which are
inherently forward-moving algorithms.

\section{Degrees of freedom and $C_{p}$ estimates.}\label{sec:sec4}
Figures 1 and 3 show all possible Lasso, Stagewise or LARS estimates
of the vector $\bbeta$ for the diabetes data. The scientists want just
a single $\widehat {\bbeta}$ of course, so we need some rule for
selecting among the possibilities. This section concerns a $C_p$-type
selection criterion, especially as it applies to the choice of LARS
estimate.

Let $\widehat \bmu = g({\bf y})$ represent a formula for estimating
$\bmu$ from the data vector ${\bf y}$. Here, as usual in regression
situations, we are considering the covariate vectors ${\bf x}_1, {\bf
  x}_2, \ldots, {\bf x}_m$ fixed at their observed values. We assume
that given the ${\bf x}$'s, ${\bf y}$ is generated according to an
homoskedastic model 
\begin{equation}\label{eqn4.1}
{\bf y} \sim (\bmu, \sigma^2 {\bf I}),
\end{equation}
meaning that the components $y_i$ are uncorrelated, with mean
$\mu_i$ and variance $\sigma^2$. Taking expectations in the identity 
\begin{equation}\label{eqn4.2}
(\widehat \mu_i - \mu_i)^2 = (y_i - \widehat \mu_i)^2 - (y_i -
\mu_i)^2 + 2 (\widehat \mu_i - \mu_i) (y_i - \mu_i),
\end{equation}
and summing over $i$, yields
\begin{equation}\label{eqn4.3}
E\biggl\{ \frac{\Vert \widehat \bmu - \bmu \Vert^2}{\sigma^2} \biggr\} =
E \biggl\{ \frac{\Vert {\bf y}-\widehat \bmu \Vert^2}{\sigma^2} - n
\biggr\} + 2 \sum_{i=1}^n \frac{\cov (\widehat \mu_i, y_i)}{\sigma^2}.
\end{equation}

The last term of (\ref{eqn4.3}) leads to a convenient definition of the {\it
  degrees of freedom} for an estimator $\widehat \bmu = g({\bf y})$, 
\begin{equation}\label{eqn4.4}
df_{\mu, \sigma^2} = \sum_{i=1}^n \cov (\widehat \mu_i, y_i) /
\sigma^2,
\end{equation}
and a $C_p$-type risk estimation formula,
\begin{equation}\label{eqn4.5}
C_p(\widehat \bmu) = \frac{\Vert {\bf y}-\widehat \bmu
  \Vert^2}{\sigma^2} - n + 2 df_{\mu, \sigma^2}.
\end{equation}
If $\sigma^2$ and $df_{\mu, \sigma^2}$ are known, $C_p(\widehat
\bmu)$ is an unbiased estimator of the true risk $E \{\Vert \widehat
\bmu - \bmu \Vert^2 / \sigma^2\}$. For linear estimators $\widehat
\bmu = M {\bf y}$, model (\ref{eqn4.1}) makes $df_{\mu, \sigma^2} = {\rm trace}
(M)$, equaling the usual definition of degrees of freedom for OLS, and
coinciding with the proposal of \citet{Ma73}.
Section~\ref{sec:sec6} of \citet{ET97} 
and Section~\ref{sec:sec7} of \citet{efro:1986}
discuss formulas (\ref{eqn4.4}) and (\ref{eqn4.5}) and their 
role in $C_p$, Akaike 
information criterion  (AIC) 
and Stein's unbiased risk 
estimated (SURE)
estimation theory, a more recent reference being \citet{ye93:_in}.

Practical use of $C_p$ formula (\ref{eqn4.5}) requires preliminary estimates of
$\bmu, \sigma^2$ and $df_{\mu, \sigma^2}$. In the numerical results
below, the usual OLS estimates $\bar {\bmu}$ and $\bar \sigma^2$ from
the full OLS model were used to calculate bootstrap estimates of
$df_{\mu, \sigma^2}$; bootstrap samples ${\bf y}^\ast$ and
replications $\widehat \bmu^\ast$ were then generated according to
\begin{equation}\label{eqn4.6}
{\bf y}^\ast \sim N(\bar {\bmu}, \bar \sigma^2) \quad \mbox{and}
\quad \widehat \bmu^\ast = g({\bf y}^\ast).
\end{equation}
Independently repeating (\ref{eqn4.6}) say $B$ times gives straightforward
estimates for the covariances in (\ref{eqn4.4}),
\begin{equation}\label{eqn4.7}
\qquad \widehat \cov_i = \frac{\sum_{b=1}^B \widehat \bmu_i^\ast(b)[\,{\bf
y}_i^\ast(b) - {\bf y}_i^\ast(\cdot)]}{B-1}, \qquad \mbox{where } 
{\bf y}^\ast (\cdot) = \frac{\sum_{b=1}^B {\bf y}^\ast(b)}{B},
\end{equation}
and then
\begin{equation}
\widehat{df} = \sum_{i=1}^n \widehat \cov_i / \bar{\sigma}^2.
\label{eqn4.8}
\end{equation}
Normality is not crucial in (\ref{eqn4.6}). Nearly the same results were
obtained using ${\bf y}^\ast = \bar {\bmu}^\ast + {\bf e}^\ast$,
where the components of ${\bf e}^\ast$ were resampled from ${\bf e} =
{\bf y} - \bar{\bmu}$.
\begin{figure}

\includegraphics{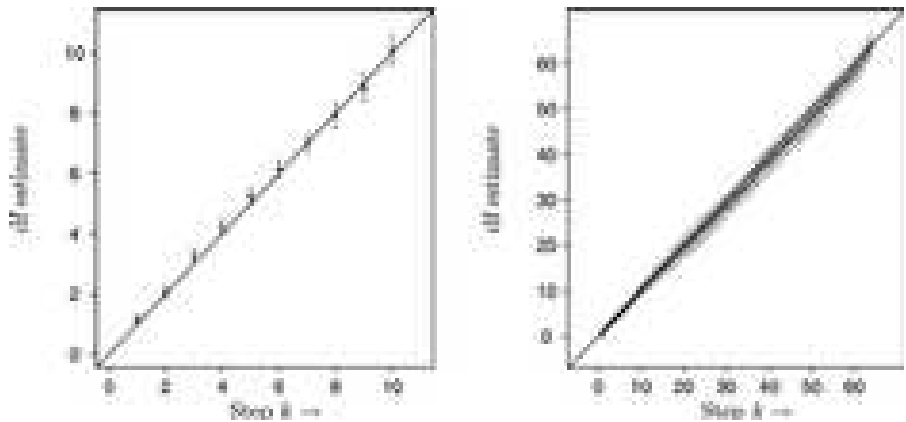}

\caption{Degrees of freedom for
    LARS estimates $\widehat {\bmu}_k$: \textup{(}{\rm left}\textup{)}  
diabetes
    study, Table $1$, $k = 1,\allowbreak  2, \ldots, m = 10$; 
{\rm (right)}
    quadratic model~\textup{(\ref{eq:insert3.1})} 
for the diabetes data, $m = 64$. Solid
    line is simple approximation $df_k = k$. Dashed lines are
    approximate $95\%$ confidence intervals for the bootstrap estimates.
    Each panel based on $B = 500$ bootstrap replications.}
\label{fig5}
\end{figure}

The left panel of Figure~\ref{fig5} shows $\widehat{df}_{\!k}$ for the diabetes
data LARS estimates~$\widehat{\bmu}_k, k = 1, 2, \ldots, m = 10$. It
portrays a startlingly simple situation that we will call the ``simple
approximation,''
\begin{equation}\label{eqn4.9}
df(\widehat{\bmu}_k) \ \dot{=} \ k.
\end{equation}
The right panel also applies to the diabetes data, but this time
with the quadratic model (\ref{eq:insert3.1}), having $m = 64$
predictors.
 We see that the simple
approximation~(\ref{eqn4.9}) is again accurate within the limits of the
bootstrap computation (\ref{eqn4.8}), where $B = 500$ replications were divided
into 10 groups of 50 each in order to calculate Student-$t$ confidence
intervals.

If (\ref{eqn4.9}) can be believed, and we will offer some evidence in its
behalf, we can estimate the risk of a $k$-step LARS estimator
$\widehat{\bmu}_k$ by
\begin{equation}\label{eqn4.11}
C_p(\widehat {\bmu}_k) \ \dot{=} \ 
\Vert {\bf y} - \widehat {\bmu}_k
\Vert^2 / \bar \sigma^2 - n + 2 k.
\end{equation}
The formula, which is the same as the $C_p$ estimate of risk for
an OLS estimator based on a subset of $k$ preselected predictor
vectors, has the great advantage of not requiring any further
calculations beyond those for the original LARS estimates. The formula
applies only to LARS, and not to Lasso or Stagewise.

\begin{figure}

\includegraphics{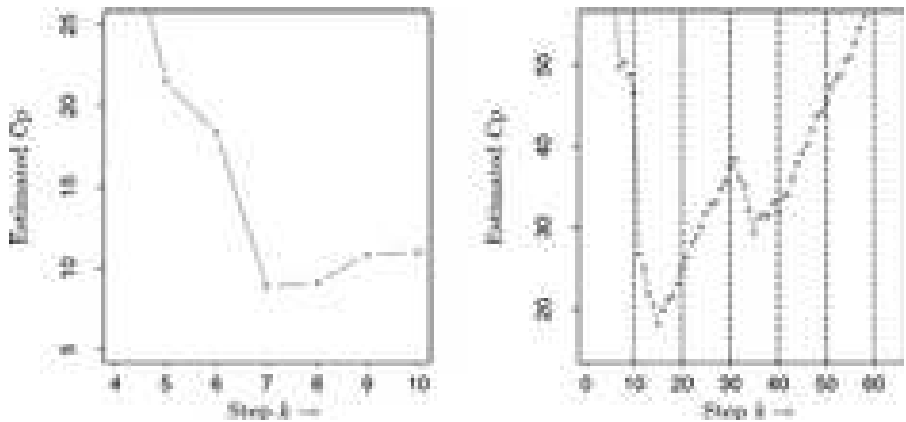}

\caption{$C_p$ estimates of
    risk \textup{(\ref{eqn4.11})} 
for the two situations of Figure~\textup{\ref{fig5}}: 
{\rm (left)}
$m = 10$ model has smallest $C_p$ at $k = 7$; {\rm (right)}  
$m = 64$ model has smallest $C_p$ at $k = 16$.}
\label{fig6}
\end{figure}

Figure~\ref{fig6} displays $C_p(\widehat \bmu_k)$ as a function of $k$ for the
two situations of Figure~\ref{fig5}. Minimum $C_p$ was achieved at steps $k =
7$ and $k = 16$, respectively. Both of the minimum $C_p$ models looked
sensible, their first several selections of ``important'' covariates
agreeing with an earlier model based on a detailed inspection of the
data assisted by medical expertise.

The simple approximation becomes a theorem in two cases.

\begin{thm}
If the covariate vectors ${\bf x}_1, {\bf x}_2, \ldots, {\bf x}_m$ are mutually
orthogonal, then the $k$-step LARS estimate $\hat \bmu_k$ has $ df(
\hat \bmu_k ) = k$. 
\end{thm}

To state the second more general setting we introduce the following 
condition. 
\begin{posit*}
For all possible subsets $X_{{\mathcal A}}$ of the full design 
matrix $X$,
\begin{equation}
\label{eq:pcc}
G_{{\mathcal A}}^{-1} \bolds{1}_{{\mathcal A}} > 0,
\end{equation}
where the inequality is taken element-wise.
\end{posit*}

The positive cone condition holds if $X$ is orthogonal. It is strictly
more general than orthogonality, but counterexamples (such as the
diabetes data) show that not all design matrices $X$ satisfy it.

It is also easy to show that LARS, Lasso and Stagewise all coincide
under the positive cone condition, so the degrees-of-freedom formula
applies to them too in this case.

\begin{thm}
\label{thm:dfcalc}
Under the positive cone condition, $ df( \hat
\bmu_k) = k$.
\end{thm}

The proof, which appears later in this section, is an application of
Stein's unbiased risk estimate (SURE) [\citet{St81}].
Suppose that $g\dvtx \mathbb{R}^n \rightarrow
\mathbb{R}^n$ is almost differentiable (see Remark A.1 in the Appendix)
and set 
$\nabla \cdot g = \sum_{i=1}^n \partial g_i / \partial x_i.$
If~${\bf y} \sim N_n( \bmu, \sigma^2 {\bf I})$, then Stein's formula states that
\begin{equation}
\label{eq:stein}
  \sum_{i=1}^n \mbox{cov} (g_i, y_i)/ \sigma^2 = 
E [\nabla \cdot g ({\bf y})].
\end{equation}
The left-hand side is $df(g)$ for the general estimator $g({\bf y})$. Focusing
specifically on LARS, it will turn out that $\nabla \cdot \hat
\bmu_k ({\bf y}) = k$ in \textit{all} situations with probability 1, but
that the continuity assumptions underlying (\ref{eq:stein})
and SURE can fail in certain
nonorthogonal cases where the positive cone condition does not hold.

 A range of simulations suggested that
the simple approximation is quite accurate even when the ${\bf x}_j$'s
are highly correlated and that it requires concerted effort at
pathology to make $df(\widehat \bmu_k)$ much different than $k$.

Stein's formula assumes normality, ${\bf y} \sim N(\bmu, \sigma^2
{\bf I})$. A cruder ``delta method'' rationale for the simple approximation
requires only homoskedasticity, 
(\ref{eqn4.1}). The geometry of Figure~\ref{fig4} implies
\begin{equation}\label{eqn4.12}
\widehat \bmu_k = \bar {\bf y}_k - {\rm cot}_k \cdot \Vert
\bar{\bf y}_{k+1} - \bar{\bf y}_k \Vert,
\end{equation}
where ${\rm cot}_k$ is the cotangent of the angle between ${\bf
 u}_k$ and ${\bf u}_{k+1}$, 
\begin{equation}\label{eqn4.13}
{\rm cot}_k = \frac{{\bf u}_k^\prime {\bf u}_{k+1}}{[1-({\bf
u}_k^\prime {\bf u}_{k+1})^2]^{1/2}}.
\end{equation}
Let ${\bf v}_k$ be the unit vector orthogonal to ${\mathcal L}(X_b)$,
the linear space spanned by the first $k$ covariates selected by LARS,
and pointing into ${\mathcal L}(X_{k+1})$ along the direction of $\bar{\bf
  y}_{k+1} - \bar{\bf y}_k$. For ${\bf y}^\ast$ near ${\bf y}$ we can
reexpress (\ref{eqn4.12}) as a locally linear transformation,
\begin{equation}\label{eqn4.14}
\widehat \bmu_k^\ast = \widehat \bmu_k + M_k ({\bf y}^\ast - {\bf y})
\qquad \mbox{with } M_k = P_k - {\rm cot}_k \cdot {\bf u}_k {\bf
 v}_k^\prime,
\end{equation}
$P_k$ being the usual projection matrix from $\mathbb{R}^n$ 
into ${\mathcal
  L}(X_k)$; (\ref{eqn4.14}) holds within a neighborhood of ${\bf y}$ such that
the LARS choices ${\mathcal L}(X_k)$ and ${\bf v}_k$ remain the same.

The matrix $M_k$ has trace$(M_k) = k$. Since the trace equals the
degrees of freedom for linear estimators, the simple approximation
(\ref{eqn4.9}) is seen to be a delta method approximation to the bootstrap
estimates (\ref{eqn4.6}) and (\ref{eqn4.7}).

It is clear that 
(\ref{eqn4.9}) $df(\widehat{\bmu}_k) \ \dot = \ k$ cannot hold for the
Lasso, since the degrees of freedom is $m$  for the full model but
the total number of steps taken can exceed $m$.
However, we have found empirically that an intuitively plausible result holds:
the degrees of freedom is well
 approximated by the number of nonzero predictors in the
model. Specifically, 
starting at step 0, let $\ell(k)$ be the index of the last model in the Lasso sequence containing
$k$ predictors. Then $df(\widehat{\bmu}_{\ell(k)})  \ \dot{=} \ k$. 
We do not yet have any mathematical support for this claim.

\subsection{Orthogonal designs.}
\label{sec:orthogonal-designs}

In the orthogonal case, we assume that ${\bf x}_j = {\bf e}_j$ 
for~$j = 1,\ldots,m$. 
The LARS algorithm then has a particularly simple form, reducing to
soft thresholding at the order statistics of the data.

To be specific, define the soft thresholding operation on a scalar
$y_1$ at threshold~$t$ by
\[
\eta(y_1; t) = 
\cases{
y_1 - t,  &\quad  \mbox{if }$y_1 > t$, \cr
0,     & \quad \mbox{if }$|y_1| \leq t$, \cr
y_1 + t,  &\quad  \mbox{if }$y_1 < -t$.}
\]
The order statistics of the absolute values of the data are denoted by
\begin{equation}
\label{eq:orderstat}
|y|_{(1)} \geq |y|_{(2)} \geq \cdots \geq |y|_{(n)} \geq
|y|_{(n+1)}:= 0.  
\end{equation}
We note that $y_{m+1}, \ldots, y_n$ do not enter into the
estimation procedure, and so we may as well assume that $m = n.$

\begin{lemma}
  \label{prop:orthog}
For an orthogonal design with ${\bf x}_j = {\bf e}_j, j = 1,\ldots, n$, 
the $k$th
LARS estimate ($0 \leq k \leq n$) is given by
\begin{eqnarray}
\label{eq:threshform}
\hat \mu_{k,i}({\bf y}) &=& 
\cases{
y_i - |y|_{(k+1)},  & \quad  \mbox{if }$y_i > |y|_{(k+1)}$, \cr
     0,                  & \quad \mbox{if }$|y_i| \leq |y|_{(k+1)}$, \cr
y_i + |y|_{(k+1)},  & \quad  \mbox{if }$y_i < - |y|_{(k+1)}$,}
\\
&=& \eta\bigl(y_i; |y|_{(k+1)} \bigr).
\end{eqnarray}
\end{lemma}

\begin{pf}
  The proof is by induction, stepping through the LARS sequence. First
  note that the LARS parameters take a simple form in the orthogonal
  setting:
\[
    G_{\mathcal A} = I_{\mathcal A}, \qquad
    A_{\mathcal A} = | {\mathcal A} |^{-1/2}, \qquad
    {\bf u}_{\mathcal A} = |{\mathcal A} |^{-1/2} \bolds{1}_{\mathcal A}, \qquad
    a_{k,j} = 0, \qquad j \notin {\mathcal A}_k.
\]
We assume for the moment that there are no ties in the order
statistics (\ref{eq:orderstat}), so that the variables enter one at a time.
Let $j(l)$ be the index corresponding to the $l$th order statistic,
$|y|_{(l)} = s_l y_{j(l)}$: we will see that ${\mathcal A}_k = \{ j(1),
\ldots, j(k) \}.$ 

We have ${\bf x}_j'{\bf y} = y_j$, and so at the first step 
LARS picks variable
$j(1)$ and sets $\hat C_1 = |y|_{(1)}.$ It is easily seen that
\[
  \hat \gamma_1 = \min_{j \neq j(1) } \bigl\{ |y|_{(1)} - |y_j| \bigr\} 
  = |y|_{(1)} - |y|_{(2)}
\]
and so
\[
  \hat \bmu_1 = \bigl[ |y|_{(1)} - |y|_{(2)} \bigr] {\bf e}_{j(1)},
\]
which is precisely (\ref{eq:threshform}) for $k=1$.

Suppose now that step $k-1$ has been completed, so that ${\mathcal A}_k = 
\{ j(1),\ldots, j(k) \}$  and (\ref{eq:threshform}) holds for $\hat
\bmu_{k-1}$. The current correlations $\hat C_k = |y|_{(k)}$ and
$\hat c_{k,j} = y_j$ for~$j \notin {\mathcal A}_k$. Since $A_k - a_{k,j} =
k^{-1/2}$, we have 
\[
  \hat \gamma_k = \min_{j \notin {\mathcal A}_k } k^{1/2} \bigl\{ |y|_{(k)} -
  |y_j| \bigr\}
\]
and
\[
\hat \gamma_k {\bf u}_k = \bigl[ |y|_{(k)} - 
|y|_{(k+1)} \bigr] \bolds{1} \{ j \in {\mathcal A}_k \}.
\]
Adding this term to $\hat \bmu_{k-1}$ yields (\ref{eq:threshform}) for
step $k$.

The argument clearly extends to the case in which there are ties in
the order statistics (\ref{eq:orderstat}): if $|y|_{(k+1)} = \cdots =
|y|_{(k+r)}$, then ${\mathcal A}_k({\bf y})$ expands by $r$ variables at step
$k+1$ and $\hat \bmu_{k+\nu}({\bf y}), \nu = 1, \ldots, r$, 
are all determined at
the same time and are equal to $\hat \bmu_{k+1}({\bf y})$. 
\end{pf}

\begin{pf*}{Proof of Theorem \textup{\ref{thm:dfcalc}} 
(\normalfont{Orthogonal case})}
The argument is particularly simple in this setting, and so worth
giving separately. First we note from (\ref{eq:threshform}) that $\hat
\bmu_k$ is continuous and Lipschitz$(1)$ and so certainly almost
differentiable. Hence~(\ref{eq:stein}) shows that we simply have
to calculate $\nabla \cdot \hat \bmu_k$. Inspection of
(\ref{eq:threshform}) shows that
\begin{eqnarray*}
\nabla \cdot \hat \bmu_k  &=& \sum_i \frac{\partial \hat
\mu_{k,i}}{\partial y_i} ({\bf y}) \\
&=& \sum_i I \bigl\{ |y_i| > |y|_{(k+1)} \bigr\} = k
\end{eqnarray*}
almost surely, that is, except for ties. This completes the proof.
\end{pf*}

\subsection{The divergence formula.}
\label{sec:divergence-formula}

While for the most general design matrices $X$, it can happen that
  $\hat \bmu_k$  fails to be almost differentiable, we will see that
  the divergence formula
\begin{equation}
    \label{eq:divergence-formula}
\nabla \cdot \hat \bmu_k ({\bf y}) = k
  \end{equation}
does hold almost everywhere. Indeed, certain authors
[e.g., \citet{meyer00}] 
have argued that the divergence $\nabla \cdot
\widehat{\bmu}$ of an estimator provides itself a useful measure of the
effective dimension of a model.

Turning to LARS, we shall say that $\hat {\bmu} ({\bf y})$ is locally linear at
a data point $y_0$ if there is some small open neighborhood of $y_0$
on which $\hat{\bmu} ({\bf y}) = M {\bf y}$ is exactly linear. Of course, the matrix
$M = M(y_0)$ can depend on $y_0$---in the case of LARS, it will be
seen to be constant on the interior of polygonal regions, with jumps across the
boundaries.  We say that a set $G$
has full measure if its complement has Lebesgue measure zero.

\begin{lemma}  \label{prop:divergence}
There is an open set $G_k$ of full measure such that, 
at all ${\bf y} \in
G_k$,  $\hat {\bmu}_k ({\bf y})$ is
locally linear  and  $\nabla \cdot \hat \bmu_k ({\bf y}) = k$.
\end{lemma}

\begin{pf}
We give here only the part of the proof that relates to actual
calculation of the divergence in (\ref{eq:divergence-formula}). The
arguments establishing continuity and local linearity are delayed to
the Appendix.

So, let us fix a point ${\bf y}$ in the interior of $G_k$. 
From Lemma
  \ref{lem-gk} in the  Appendix, this means that
  near ${\bf y}$ the active set ${\mathcal A}_k({\bf y})$ is locally constant, that a
  single variable enters at the next step, this variable being the
  same near ${\bf y}$. In addition, $\hat\bmu_k({\bf y})$ is locally linear, and
  hence in particular differentiable. Since $G_k \subset G_l$ for $l <
  k$, the same story applies at all previous steps and we have
  \begin{equation}
    \label{eq:muk1}
\hat \bmu_k ({\bf y}) = \sum_{l=1}^k \gamma_l({\bf y}) {\bf u}_l.
  \end{equation}
Differentiating the $j$th component of vector $\hat \bmu_k({\bf y})$ yields 
\[
  \frac{\partial \hat \mu_{k,j}}{\partial y_i} ({\bf y}) =
  \sum_{l=1}^k \frac{\partial \gamma_l ({\bf y})}{\partial y_i} u_{l,j}.
\]
In particular, for the divergence
\begin{equation}
  \label{eq:divergence}
    \nabla \cdot \hat \bmu_k ({\bf y})  = 
  \sum_{i=1}^n \ \frac{\partial \hat \mu_{k,i}}{\partial y_i} 
  = \sum_{l=1}^k \ \langle \nabla \gamma_l, {\bf u}_l \rangle,
\end{equation}
the brackets indicating inner product.

The active set is ${\mathcal A}_k = \{ 1,2,\ldots, k \}$ 
and ${\bf x}_{k+1}$ is
the variable to enter next.  
For~$k \geq 2,$ write $\bdelta_k = {\bf x}_l -
{\bf x}_k$ for any choice $l < k$---as remarked in the Conventions
 in the Appendix, the choice of $l$ is immaterial (e.g., $l=1$ for
definiteness).  Let~$b_{k+1} = \langle \bdelta_{k+1}, 
{\bf u}_k \rangle$,
which is nonzero, as argued in the proof of Lemma~\ref{lem-gk}.
As shown in (\ref{eq:gammaky}) in the
Appendix, (2.13) can be rewritten
\begin{equation}
\label{eq:gammaky1}
\gamma_k ({\bf y}) = b_{k+1}^{-1} \langle \bdelta_{k+1}, {\bf y} 
- \hat \bmu_{k-1}
\rangle. 
\end{equation}
For $k \geq 2$, define the linear space of vectors equiangular with
the active set
\[
  {\mathcal L}_k = {\mathcal L}_k({\bf y}) = 
  \bigl\{{\bf u} \dvtx \langle {\bf x}_1, {\bf u} \rangle = \cdots = \langle
  {\bf x}_k, {\bf u} \rangle \ 
  \mbox{ for } {\bf x}_l \mbox{ with } l \in {\mathcal A}_k({\bf y}) 
\bigr\}.
\]
[We may drop the dependence on ${\bf y}$ since ${\mathcal A}_k({\bf y})$ is locally fixed.]
Clearly $\mbox{dim} \, {\mathcal L}_k = n - k + 1$ and 
\begin{equation}
  \label{eq:containment}
  {\bf u}_k \in {\mathcal L}_k, \qquad {\mathcal L}_{k+1} 
\subset {\mathcal L}_k.
\end{equation}

We shall now verify that, for each $k \geq 1$,
\begin{equation}
  \label{eq:claimk}
  \langle \nabla \gamma_k, {\bf u}_k \rangle = 1 \quad
  \mbox{and} \quad \langle \nabla \gamma_k, {\bf u} \rangle = 0 \qquad 
  \mbox{for }  {\bf u} \in {\mathcal L}_{k+1}.
\end{equation}
Formula (\ref{eq:divergence}) shows that this suffices to prove 
Lemma~\ref{prop:divergence}. 

First, for $k = 1$ we have $\gamma_1({\bf y}) = b_2^{-1} \langle \bdelta_2, {\bf y}
\rangle$ and $\langle \nabla \gamma_1, {\bf u} \rangle = b_2^{-1} \langle
\bdelta_2, {\bf u} \rangle$, and that
\[
  \langle \bdelta_2, {\bf u} \rangle 
= \langle {\bf x}_1 - {\bf x}_2, {\bf u} \rangle = 
\cases{
    b_2,  &\quad  \mbox{if }  ${\bf u} = {\bf u}_1$, \cr
    0,   & \quad \mbox{if }  ${\bf u} \in {\mathcal L}_2$. }      
\]
Now, for general $k$, combine (\ref{eq:gammaky1}) and (\ref{eq:muk1}):
\[
b_{k+1} \gamma_k({\bf y}) = \langle \bdelta_{k+1},{\bf y} \rangle - 
                \sum_{l=1}^{k-1} \langle \bdelta_{k+1},
                {\bf u}_l \rangle \gamma_l ({\bf y}),
\]
and hence 
\[
  b_{k+1} \langle \nabla \gamma_k, {\bf u} \rangle = 
  \langle \bdelta_{k+1}, {\bf u} \rangle -
  \sum_{l=1}^{k-1} \langle \bdelta_{k+1}, {\bf u}_l \rangle  
  \langle \nabla \gamma_l, {\bf u} \rangle.
\]
From the definitions of $b_{k+1}$ and ${\mathcal L}_{k+1}$ we have
\[
  \langle \bdelta_{k+1}, {\bf u} \rangle = 
  \langle {\bf x}_l - {\bf x}_{k+1} \rangle = 
  \cases{
    b_{k+1},  &\quad  \mbox{if } ${\bf u} = {\bf u}_k$, \cr
    0,   & \quad \mbox{if }  ${\bf u} \in {\mathcal L}_{k+1}$.}
\]
Hence the truth of (\ref{eq:claimk}) for step $k$ follows from its
truth at step $k-1$ because of the containment properties
(\ref{eq:containment}). 
\end{pf}

\subsection{Proof of Theorem \textup{\ref{thm:dfcalc}}.}
\label{sec:proof-theorem-}

To complete the proof of Theorem \ref{thm:dfcalc}, we state the
following regularity result, proved in the Appendix.

\begin{lemma} \label{prop:regularity}
Under the positive cone condition, $\hat \bmu_k ({\bf y})$ is continuous
and almost differentiable.
\end{lemma}

This guarantees that Stein's formula (\ref{eq:stein}) is
  valid for $\hat \bmu_k$ under the positive cone condition, so the
  divergence formula of Lemma
  \ref{prop:divergence} then immediately yields Theorem \ref{thm:dfcalc}.

\section{LARS and Lasso properties.}\label{sec:sec5} 
The LARS and Lasso
algorithms are described more carefully in this section, with an eye
toward fully understanding their relationship. Theorem~1 of Section~\ref{sec:sec3} will be
verified. The latter material overlaps results in
\citet{osborne00}, particularly in their Section~4.
Our point of
view here allows the Lasso to be described as a quite simple
modification of LARS, itself a variation of traditional Forward
Selection methodology, and in this sense should be more accessible to
statistical audiences. In any case we will stick to the language of
regression and correlation rather than convex optimization, though
some of the techniques are familiar from the optimization literature.

The results will be developed in a series of lemmas, eventually
lending to a proof of Theorem~1 and its generalizations. The first
three lemmas refer to attributes of the LARS procedure that are not
specific to its Lasso modification.

Using notation as in (\ref{eqn2.17})--(\ref{eqn2.20}), suppose LARS has completed step
$k-1$, giving estimate $\widehat \bmu_{k-1}$ and active set ${\mathcal
  A}_k$ for step $k$, with covariate ${\bf x}_k$ the newest addition
to the active set.

\begin{lemma}
\label{lem:lem1}
If ${\bf x}_k$ is the only addition to the
active set at the end of step~$k-1$, then the coefficient vector $w_k
= A_k {\mathcal G}_k^{-1} \bolds{1}_k$ for the equiangular vector ${\bf u}_k = X_k
w_k$,~\textup{(\ref{eqn2.6})}, has its $k$th component $w_{kk}$ agreeing in sign with
the current correlation $c_{kk} = {\bf x}_k^\prime({\bf y} - \widehat
\bmu_{k-1})$. Moreover, the regression vector $\widehat \bbeta_k$ for
$\widehat \bmu_k = X \widehat \bbeta_k$ has its $k$th component
$\widehat \beta_{kk}$ agreeing in sign with $c_{kk}$.
\end{lemma}

Lemma~\ref{lem:lem1} says that new variables {\it enter} the LARS active set in the
``correct'' direction, a weakened version of the Lasso requirement
(\ref{eqn3.1}). This will turn out to be a crucial connection for the
LARS--Lasso relationship.

\begin{pf*}{Proof of Lemma~4}
 The case $k=1$ is apparent. Note that since
\[
X_k^\prime({\bf y} - \widehat \bmu_{k-1})=\widehat C_k{\bf 1}_k,
\]
(\ref{eqn2.20}), from (\ref{eqn2.6}) we have
\begin{equation}
  \label{eqn.trev1}
  w_k=A_k\widehat C_k^{-1}[(X_k^\prime X_k)^{-1}X_k^\prime({\bf
  y}-\widehat \bmu_{k-1})]:=A_k\widehat C_k^{-1} w_k^*.
\end{equation}
The term in square braces is the least squares coefficient vector in the
regression of the current residual on $X_k$, and the term preceding it
is positive.

Note also that 
\begin{equation}
  \label{eqn.trev2}
  X_k^\prime ({\bf y}-\bar{\bf y}_{k-1})=({\bf 0},\delta)^\prime\qquad
  \mbox{with }\delta > 0,
\end{equation}
since $X_{k-1}^\prime({\bf y}-\bar{\bf y}_{k-1})=\bf{0}$ by definition
(this $\bf 0$ has $k-1$ elements), and $c_k(\gamma)={\bf x}_k^\prime
({\bf y}-\gamma{\bf u}_{k-1})$ decreases more slowly in $\gamma$ than
$c_j(\gamma)$ for $j\in {\mathcal A}_{k-1}$:
\begin{equation}\label{eqn.trev3}
c_k(\gamma)\cases{
< c_j(\gamma), & \quad \mbox{for }$\gamma <
\widehat\gamma_{k-1}$,\cr
= c_j(\gamma)=\widehat C_k,& \quad  \mbox{for }$\gamma =
\widehat\gamma_{k-1}$,\cr
> c_j(\gamma), & \quad \mbox{for }$\widehat\gamma_{k-1}<\gamma <
\bar\gamma_{k-1}$.}
\end{equation}

Thus
\begin{eqnarray}
\widehat w_k^* &=&  (X_k^\prime X_k)^{-1}X_k^\prime({\bf
  y}-  \bar{\bf y}_{k-1} + \bar{\bf y}_{k-1} -\widehat \bmu_{k-1})
  \label{eqn.trev4a}\\
&=& (X_k^\prime X_k)^{-1}
\left(
  \begin{array}{c}
\bf 0\\ \delta
  \end{array}
\right)
 +(X_k^\prime
X_k)^{-1}X_k^\prime[(\bar{\gamma}_{k-1}-\widehat\gamma_{k-1}){\bf
 u}_{k-1}].
\label{eqn.trev4b}
\end{eqnarray}
The $k$th element of $\widehat w^*_k$ is positive, because it is in the first
 term in (\ref{eqn.trev4b}) [$(X_k^\prime X_k)$ is positive definite],
 and in the second term it is $0$ 
since ${\bf u}_{k-1}\in {\mathcal L}(X_{k-1})$.
 
 This proves the first statement in Lemma~\ref{lem:lem1}. The second
follows from
\begin{equation}
\widehat \beta_{kk} = \widehat \beta_{k-1,k} + \widehat \gamma_k
w_{kk},
\end{equation}
and $\widehat \beta_{k-1,k} = 0, \ {\bf x}_k$ not being active
before step $k$.
\end{pf*}

Our second lemma interprets the quantity 
$A_{\mathcal A} = ({\bf 1}^\prime {\mathcal
  G}_{\mathcal A}^{-1} {\bf 1})^{-1/2}$, (\ref{eqn2.4}) and (\ref{eqn2.5}). 
Let~${\mathcal
  S}_{\mathcal A}$ indicate the extended simplex generated by the columns
of $X_{\mathcal A}$,
\begin{equation}
{\mathcal S}_{\mathcal A} = \Biggl\{ {\bf v} = \sum_{j \in {\mathcal A}} s_j
{\bf x}_j P_j \dvtx \sum_{j \in {\mathcal A}} P_j = 1 \Biggr\},
\end{equation}
``extended'' meaning that the coefficients $P_j$ are allowed to
be negative. 

\begin{lemma}
\label{lem:lem2}
The point in ${\mathcal S}_{\mathcal A}$ nearest the
origin is
\begin{equation}
\label{eqn5.12}
{\bf v}_{\mathcal A} = A_{\mathcal A} {\bf u}_{\mathcal A} = 
A_{\mathcal A} X_{\mathcal A}
w_{\mathcal A} \qquad \mbox{where } w_{\mathcal A} = 
A_{\mathcal A} {\mathcal G}_{\mathcal A}^{-1}
{\bf 1}_{\mathcal A},
\end{equation}
with length $\Vert {\bf v}_{\mathcal A} \Vert = A_{\mathcal A}$. If
${\mathcal A} \subseteq {\mathcal B}$, then $A_{\mathcal A} \geq A_{\mathcal B}$, the
largest possible value being $A_{\mathcal A} = 1$ for ${\mathcal A}$ a
singleton. 
\end{lemma}

\begin{pf}
For any ${\bf v}\in {\mathcal S}_{\mathcal A}$, the squared distance to the 
origin is
$\Vert X_{\mathcal A}P\Vert^2=P^\prime{\mathcal G}_{\mathcal A}P$.
Introducing a Lagrange multiplier to enforce the summation constraint,
we differentiate
\begin{equation}
  \label{eq:th21}
  P^\prime{\mathcal G}_{\mathcal A}P -\lambda({\bf 1}_{\mathcal A}^\prime P - 1),
\end{equation}
and find that the minimizing  
$P_{\mathcal A}=\lambda {\mathcal G}_{\mathcal A}^{-1} {\bf 1}_{\mathcal A}$. 
Summing, 
we get $\lambda{\bf 1}_{\mathcal A}^\prime
{\mathcal G}_{\mathcal A}^{-1}
{\bf 1}_{\mathcal A}=1$, and hence 
\begin{equation}\label{eqn5.14}
P_{\mathcal A} = A_{\mathcal A}^2  {\mathcal G}_{\mathcal A}^{-1} {\bf 
1}_{\mathcal A}
= A_{\mathcal A}  w_{\mathcal A}.
\end{equation}
Hence  ${\bf v}_{\mathcal A} = X_{\mathcal A} P_{\mathcal A}
\in {\mathcal S}_{\mathcal A}$ and
\begin{equation}
\Vert {\bf v}_{\mathcal A} \Vert^2 = 
P_{\mathcal A}^\prime {\mathcal G}_{\mathcal
  A}^{-1} P_{\mathcal A} = A_{\mathcal A}^4 
{\bf 1}_{\mathcal A}^\prime {\mathcal G}_{\mathcal
  A}^{-1} {\bf 1}_{\mathcal A} = A_{\mathcal A}^2,
\end{equation}
verifying (\ref{eqn5.12}). If ${\mathcal A} \subseteq {\mathcal B}$, 
then ${\mathcal
  S}_{\mathcal A} \subseteq {\mathcal S}_{\mathcal B}$, so the nearest distance
$A_{\mathcal B}$ must be equal to or less than the nearest distance 
$A_{\mathcal
  A}$. $A_{\mathcal A}$ obviously equals $1$ if and only if~${\mathcal A}$ has
only one member.
\end{pf}


The LARS algorithm and its various modifications proceed in piecewise
linear steps. For $m$-vectors $\widehat \bbeta$ and {\bf d}, let
\begin{equation}\label{eqn5.16}
\bbeta(\gamma) = \widehat \bbeta + \gamma {\bf d} \quad \mbox{and}
\quad S (\gamma) = \Vert {\bf y} - X \bbeta (\gamma) \Vert^2.
\end{equation}

\begin{lemma}
  \label{lem:lem3}
Letting $\widehat {\bf c} = X^\prime ({\bf y}
- X \widehat \bbeta)$ be the current correlation vector at $\widehat
\bmu = X \widehat \bbeta$, 
\begin{equation}
\label{eqn5.17}
S (\gamma) - S(0) = -2\, \widehat {\bf c}^{\,\prime} {\bf d} \gamma + {\bf
  d}^\prime X^\prime X {\bf d} \gamma^2.
\end{equation}
\end{lemma}

\begin{pf}
${\mathcal S}(\gamma)$ is a quadratic function of
$\gamma$, with first two derivatives 
at $\gamma = 0$,
\begin{equation}\label{eqn5.18}
\dot S(0) = -2\, \widehat {\bf c\,}^\prime {\bf d} \quad \mbox{and}
\quad \ddot S (0) = 2 {\bf d}^\prime X^\prime X {\bf d}.
\end{equation}
\upqed\end{pf}

The remainder of this section concerns the LARS--Lasso relationship.
Now $\widehat \bbeta = \widehat \bbeta (t)$ will indicate a Lasso
solution (\ref{eqn1.5}), and likewise $\widehat \bmu = \widehat \bmu (t) =
X \widehat \bbeta(t)$. Because $S(\widehat \bbeta)$ and $T(\widehat
\bbeta)$ are both convex functions of $\widehat \bbeta$, with $S$
strictly convex, standard results show that $\widehat \bbeta(t)$ and
$\widehat \bmu(t)$ are unique and continuous functions of $t$.

For a given value of $t$ let
\begin{equation}
{\mathcal A} = \{j \dvtx \widehat \beta_j(t) \ne 0\} .
\end{equation}
We will show later that ${\mathcal A}$ is also the active set that
determines the equiangular direction ${\bf u}_{\mathcal A}$, (\ref{eqn2.6}), for
the LARS--Lasso computations.

We wish to characterize the track of the Lasso solutions $\widehat
\bbeta(t)$ or equivalently of $\widehat \bmu(t)$ as $t$ increases from
0 to its maximum effective value. Let ${\mathcal T}$ be an open interval
of the $t$ axis, with infimum $t_0$, within which the set ${\mathcal A}$
of nonzero Lasso coefficients $\widehat \beta_j(t)$ remains constant.

\begin{lemma}
\label{lem:lem4}
The Lasso estimates $\widehat \bmu (t)$
satisfy
\begin{equation}
\label{eqn5.20}
\widehat \bmu (t) = \widehat \bmu (t_0) + A_{\mathcal A} (t - t_0) {\bf
  u}_{\mathcal A}
\end{equation}
for $t \in {\mathcal T}$, where ${\bf u}_{\mathcal A}$ is the equiangular
vector $X_{\mathcal A} w_{\mathcal A}, w_{\mathcal A} 
= A_{\mathcal A} {\mathcal G}_{\mathcal
  A}^{-1} {\bf 1}_{\mathcal A}$, \textup{(\ref{eqn2.7})}.
\end{lemma}

\begin{pf}
The lemma says that, for $t$ in ${\mathcal T}$,
$\widehat \bmu(t)$ moves linearly along the equiangular vector ${\bf
  u}_{\mathcal A}$ determined by ${\mathcal A}$. We can also state this in
terms of the nonzero regression coefficients $\widehat \beta_{\mathcal
 A}(t)$,
\begin{equation}\label{eqn5.21}
\widehat \beta_{\mathcal A}(t) = \widehat \beta_{\mathcal A}(t_0) 
+ S_{\mathcal A}
A_{\mathcal A}(t - t_0) w_{\mathcal A},
\end{equation}
where $S_{\mathcal A}$ is the diagonal matrix with diagonal elements
$s_j$, $j \in {\mathcal A}$. 
[$S_{\mathcal A}$ is needed in (\ref{eqn5.21}) because
definitions (\ref{eqn2.4}), (\ref{eqn2.10}) require $\widehat \bmu(t) = X \widehat
\bbeta(t) = X_{\mathcal A} S_{\mathcal A} \widehat \beta_{\mathcal 
A}(t)$.]

Since $\widehat \bbeta(t)$ satisfies (\ref{eqn1.5}) and has nonzero 
set ${\mathcal A}$, it also minimizes 
\begin{equation}
S (\widehat \beta_{\mathcal A}) = \Vert {\bf y} 
- X_{\mathcal A} S_{\mathcal A}
\widehat \beta_{\mathcal A} \Vert^2
\end{equation}
subject to
\begin{equation}\label{eqn5.23}
\sum_{\mathcal A} s_j \widehat \beta_j = t \quad \mbox{and} \quad
\sign (\widehat \beta_j) = s_j \qquad \mbox{for } j \in
{\mathcal A}.
\end{equation}
[The inequality in (\ref{eqn1.5}) can be replaced by $T (\widehat \bbeta) =
t$ as long as $t$ is less than $\sum \vert \bar \beta_j \vert$ for
the full $m$-variable OLS solution $\bar {\bbeta}_m$.] 
Moreover, the
fact that the minimizing point $\widehat \beta_{\mathcal A}(t)$ occurs
strictly {\it inside} the simplex (\ref{eqn5.23}), combined with the strict
convexity of $S(\widehat \beta_{\mathcal A})$, implies we can drop the
second condition in (\ref{eqn5.23}) so that~$\widehat \beta_{\mathcal A}(t)$ solves
\begin{equation}\label{eqn5.24}
\mbox{minimize}\quad \{S(\widehat \beta_{\mathcal A})\} \qquad 
\mbox{subject to}
\quad \sum_{\mathcal A} s_j \widehat \beta_j = t.
\end{equation}

Introducing a Lagrange multiplier,  
(\ref{eqn5.24}) becomes
\begin{equation}
\label{eq:th1}
\mbox{minimize}\quad \tfrac{1}{2} \Vert {\bf y}-
X_{\mathcal A}S_{\mathcal
A}\widehat{\beta}_{\mathcal A}\Vert^2+
\lambda\sum_{\mathcal A}s_j\widehat\beta_j.
\end{equation}
Differentiating we get
\begin{equation}
  \label{eq:th2}
  -S_{\mathcal A}X_{\mathcal A}^\prime({\bf y} -X_{\mathcal A}S_{\mathcal
   A}\widehat\beta_{\mathcal A})+\lambda S_{\mathcal A}{\bf 1}_{\mathcal A}=0.
\end{equation}

Consider two values $t_1$ and $t_2$ in ${\mathcal T}$ with
$t_0<t_1<t_2$. Corresponding to each of these are values for the
Lagrange multiplier $\lambda$ such that
$\lambda_1>\lambda_2$, and solutions $\widehat\beta_{\mathcal A}(t_1)$ and $\widehat\beta_{\mathcal A}(t_2)$.
Inserting these into (\ref{eq:th2}), differencing and premultiplying by~$S_{\mathcal A}$
we get
\begin{equation}
  \label{eq:th3}
  X_{\mathcal A}^\prime X_{\mathcal A}S_{\mathcal A}\bigl(\widehat\beta_{\mathcal
  A}(t_2)-\widehat\beta_{\mathcal
  A}(t_1)\bigr)=(\lambda_1-\lambda_2){\bf 1}_{\mathcal A}.
\end{equation}
Hence 
\begin{equation}
\label{eq:th4}
\widehat \beta_{\mathcal A}(t_2) - \widehat \beta_{\mathcal A}(t_1) = (\lambda_1 -
\lambda_2) S_{\mathcal A} {\mathcal G}_{\mathcal A}^{-1} {\bf 1}_{\mathcal A}. 
\end{equation}
However, 
$s_{\mathcal A}^\prime [(\widehat \beta_{\mathcal A}(t_2) - \widehat
\beta_{\mathcal A}(t_1)] = t_2 - t_1$ 
according to the Lasso definition, so
\begin{equation}
\qquad \quad \hspace*{5pt}t_2-t_1 = 
(\lam_1 - \lam_2) s_{\mathcal A}^\prime S_{\mathcal A} {\mathcal
  G}_{\mathcal A}^{-1} {\bf 1}_{\mathcal A} = 
(\lam_1 - \lam_2) {\bf 1}_{\mathcal A}^\prime
  {\mathcal G}_{\mathcal A}^{-1} {\bf 1}_{\mathcal A} = 
(\lam_1 - \lam_2) A_{\mathcal A}^{-2}
\end{equation}
and
\begin{equation}
\widehat \beta_{\mathcal A}(t_2) - \widehat \beta_{\mathcal A}(t_1) = S_{\mathcal A}
A_{\mathcal A}^2(t_2 - t_1) {\mathcal G}_{\mathcal A}^{-1}{\bf 1}_{\mathcal A}
= S_{\mathcal A} A_{\mathcal A}(t - t_1) w_{\mathcal A} .
\end{equation}

Letting $t_2=t$ and $t_1 \rightarrow t_0$ gives 
(\ref{eqn5.21}) by the continuity of
$\widehat \bbeta(t)$, and finally~(\ref{eqn5.20}). Note that 
(\ref{eqn5.20}) implies
that the maximum absolute correlation $\widehat C(t)$ equals $\widehat
C(t_0) - A_{\mathcal A}^2 (t - t_0)$, so that $\widehat C(t)$ is a
piecewise linear decreasing function of the Lasso parameter $t$.
\end{pf}

The Lasso solution $\widehat \bbeta(t)$ occurs on the surface of the
diamond-shaped convex polytope
\begin{equation}
{\mathcal D} (t) = \Bigl\{\bbeta \dvtx 
\sum \vert \beta_j \vert \leq t \Bigr\}, 
\end{equation}
${\mathcal D}(t)$ increasing with $t$. Lemma~\ref{lem:lem4} says that, 
for $t \in
{\mathcal T}$, $\widehat \bbeta(t)$ moves linearly along edge ${\mathcal A}$
of the polytope, the edge having $\beta_j = 0$ for $j \notin {\mathcal
  A}$. Moreover the regression estimates $\widehat \bmu(t)$ move in
the LARS equiangular direction ${\bf u}_{\mathcal A}$, (\ref{eqn2.6}). It remains
to show that ``${\mathcal A}$'' changes according to the rules of Theorem
1, which is the purpose of the next three lemmas.

\begin{lemma}
  \label{lem:lem5}
A Lasso solution $\widehat \bbeta$ has
\begin{equation}
\label{eqn5.37}
\widehat c_j = \widehat C \cdot \sign(\widehat \beta_j) \qquad
\mbox{for }  j \in {\mathcal A},
\end{equation}
where $\widehat c_j$ equals the current correlation ${\bf
  x}_j^\prime ({\bf y} - \widehat \bmu) = {\bf x}_j^\prime ({\bf y} -
X \widehat \bbeta)$.
In particular, this implies that 
\begin{equation}
\label{eqn5.31}
\sign (\widehat \beta_j) = \sign (\widehat c_j) \qquad
\mbox{for } j \in {\mathcal A}.
\end{equation}
\end{lemma}

\begin{pf}
This follows immediately from (\ref{eq:th2}) by noting that the $j$th
element of the left-hand side is $\widehat c_j$, and the right-hand
side is $\lambda\cdot\mbox{sign}(\widehat \beta_j)$ for $j\in \mathcal{A}$. 
Likewise $\lambda=|\widehat c_j|=\widehat C$.
\end{pf}

\begin{lemma}
  \label{lem:lem6}
Within an interval ${\mathcal T}$ of constant
nonzero set ${\mathcal A}$, and also at~$t_0 = \inf({\mathcal T})$, 
the Lasso
current correlations $c_j (t) = {\bf x}_j^\prime ({\bf y} - \widehat
\bmu (t))$ satisfy
\[
\vert c_j(t) \vert = \widehat C(t) \equiv \max\{\vert c_{\ell}(t)
\vert\} \qquad \mbox{for }  j \in {\mathcal A} 
\]
and
\begin{equation}
\label{eqn5.38}
\vert c_j(t) \vert \leq \widehat C(t) \qquad \mbox{for }  j \notin 
{\mathcal A}. 
\end{equation}
\end{lemma}

\begin{figure}[b]

\includegraphics{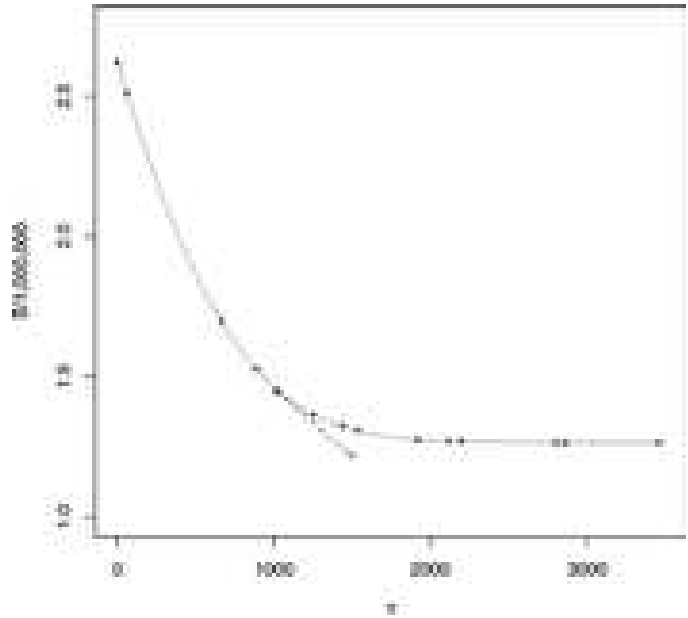}

\caption{Plot of $S$ versus $T$
    for Lasso applied to diabetes data; points indicate the $12$
    modified LARS steps of Figure~\textup{\ref{fig1}}; 
triangle is $(T, S)$ boundary
    point at $t = 1000$; dashed arrow is tangent at $t = 1000$,
    negative slope $R_t$, \textup{(\ref{eqn5.39})}. 
The $(T, S)$ curve is a decreasing,
    convex, quadratic spline.}
\label{fig7}
\end{figure}

\begin{pf}
Equation (\ref{eqn5.37}) says that the $\vert c_j(t) \vert$ have
identical values, say $\widehat C_t$, for~$j \in {\mathcal A}$. It remains
to show that $\widehat C_t$ has the extremum properties indicated 
in~(\ref{eqn5.38}). For an $m$-vector {\bf d} we define $\bbeta (\gamma) =
\widehat \bbeta(t) + \gamma {\bf d}$ and $S(\gamma)$ as in
(\ref{eqn5.16}), likewise $T(\gamma) = \sum \vert \beta_j(\gamma)
\vert$, and
\begin{equation}\label{eqn5.39}
R_t(d) = - \dot S (0) / \dot T(0).
\end{equation}
Again assuming $\widehat \beta_j > 0$ for $j \in {\mathcal A}$, by
redefinition of ${\bf x}_j$ if necessary, (\ref{eqn5.18}) 
and (\ref{eqn5.37}) yield
\begin{equation}\label{eqn5.40}
R_t({\bf d}) = 2 \Biggl[\widehat C_t \sum_{\mathcal A} d_j + \sum_{{\mathcal
    A}^c} c_j(t) d_j \Biggr] \bigg/ \Biggl[\sum_{\mathcal A} d_j +
    \sum_{{\mathcal A}^c} \vert d_j \vert \Biggr].\
\end{equation}
If $d_j = 0$ for $j \notin {\mathcal A}$, and $\sum d_j \ne 
0$,
\begin{equation}\label{eqn5.41}
R_t ({\bf d}) = 2 \widehat C_t,
\end{equation}
while if ${\bf d}$ has only component $j$ nonzero we can make 
\begin{equation}
R_t ({\bf d}) = 2 \vert c_j(t) \vert.
\end{equation}

According to Lemma~\ref{lem:lem4} the Lasso solutions for $t \in {\mathcal T}$ use
$d_{\mathcal A}$ proportional to $w_{\mathcal A}$ with $d_j = 0$ for $j 
\notin {\mathcal A}$, so
\begin{equation}
R_t \equiv R_t(w_{\mathcal A})
\end{equation}
is the downward slope of the curve ($T, S(T))$ at $T = t$, and by
the definition of the Lasso must maximize $R_t({\bf d})$. This shows
that $\widehat C_t = \widehat C(t)$, and verifies (\ref{eqn5.38}), which also
holds at $t_0 = \inf({\mathcal T})$ by the continuity of the current
correlations.
\end{pf}

We note that Lemmas \ref{lem:lem4}--\ref{lem:lem6} follow relatively
easily from the Karush--Kuhn--Tucker conditions for optimality for the
quadratic programming Lasso problem [\citet{osborne00}]; 
we have chosen a
more geometrical argument here to demonstrate the nature of the
Lasso path.

Figure~\ref{fig7} shows the $(T, S)$ curve corresponding to the Lasso estimates
in Figure~\ref{fig1}. The arrow indicates the tangent to the curve at $t =
1000$, which has downward slope $R_{1000}$. The argument above relies
on the fact that $R_t({\bf d})$ cannot be greater than $R_t$, or else
there would be $(T, S)$ values lying below the optimal curve. Using
Lemmas 3 and 4 it can be shown that the $(T, S)$ curve is always
convex, as in Figure~\ref{fig7}, being a quadratic spline with $\dot S (T) = -2
\widehat C(T)$ and $\ddot S (T) = 2 A_{\mathcal A}^2$.

We now consider in detail the choice of active set at a breakpoint of
the piecewise linear Lasso path. Let $t = t_0$ indicate 
such a point,
$t_0 = \inf ({\mathcal T})$ as in Lemma~\ref{lem:lem6}, 
with Lasso regression vector
$\widehat \bbeta$, prediction estimate $\widehat \bmu = X \widehat
\bbeta$, current correlations $\widehat {\bf c} = X^\prime ({\bf y} -
\widehat \bmu)$, $s_j = \sign (\widehat c_j)$ 
and maximum
absolute correlation $\widehat C$. Define
\begin{equation}\label{eqn5.44}
{\mathcal A}_1 = \{j \dvtx \widehat \beta_j \ne  0\}, \qquad 
{\mathcal A}_0 = \{j \dvtx
\widehat \beta_j = 0 \mbox{ and }  \vert \widehat c_j
 \vert =
\widehat C\},
\end{equation}
${\mathcal A}_{10} = {\mathcal A}_1 \cup {\mathcal A}_0$ and 
${\mathcal A}_2 =
{\mathcal A}_{10}^c$, and take $\bbeta(\gamma) = \widehat \bbeta + \gamma
{\bf d}$ for some $m$-vector {\bf d}; also $S(\gamma) = \Vert {\bf y}
- X \bbeta (\gamma) \Vert^2$ and $T(\gamma) = \sum \vert
\beta_j(\gamma)\vert$.

\begin{lemma}
\label{lem:lem7}
 The negative slope \textup{(\ref{eqn5.39})} at $t_0$ is bounded
by $2 \widehat C$, 
\begin{equation}
\label{eqn5.45}
R({\bf d}) = - \dot S (0) / \dot T (0) \leq 2 \widehat C,
\end{equation}
with equality only if $d_j = 0$ for $j \in {\mathcal A}_2$. If so,
the differences $\Delta S = S(\gamma) - S(0)$ and $\Delta T = T
(\gamma) - T(0)$ satisfy
\begin{equation}\label{eqn5.46}
\Delta S = -2 \widehat C \Delta T + L ({\bf d})^2 \cdot (\Delta T)^2,
\end{equation}
where
\begin{equation}
L ({\bf d}) = \Vert X {\bf d} / d_+ \Vert.
\label{eqn5.47}
\end{equation}
\end{lemma}

\begin{pf}
We can assume $\widehat c_j \geq 0$ for all
$j$, by redefinition if necessary, so $\widehat \beta_j \geq 0$
according to Lemma~\ref{lem:lem5}. Proceeding as in (\ref{eqn5.40}),
\begin{equation}\label{eqn5.48}
R({\bf d}) = 2 \widehat C \Biggl[\, \sum_{{\mathcal A}_{10}} d_j +
\sum_{{\mathcal A}_2} (\widehat c_j / \widehat C) d_j \Biggr] \bigg/ 
\Biggl[
\sum_{{\mathcal A}_1} d_j + \sum_{{\mathcal A}_0 \cup {\mathcal A}_2} \vert d_j
\vert \Biggr].
\end{equation}
We need $d_j \geq 0$ for $j \in {\mathcal A}_0 \cup {\mathcal A}_2$ in
order to maximize (\ref{eqn5.48}), in which case
\begin{equation}
R ({\bf d}) = 2 \widehat C \Biggl[ \,\sum_{{\mathcal A}_{10}} d_j +
\sum_{{\mathcal A}_2} (\widehat c_j / \widehat C) d_j \Biggr] \bigg/ 
\Biggl[
\,\sum_{{\mathcal A}_{10}} d_j + \sum_{{\mathcal A}_2} d_j \Biggr].
\end{equation}
This is $< 2 \widehat C$ unless $d_j = 0$ for $j \in {\mathcal A}_2$,
verifying (\ref{eqn5.45}), and also implying
\begin{equation}
T(\gamma) = T(0) + \gamma \sum_{{\mathcal A}_{10}} d_j.
\end{equation}
The first term on the right-hand side of (\ref{eqn5.17}) is then $-2 \widehat C
(\Delta T)$, while the second term equals $({\bf d} / d_+)^\prime
X^\prime X ({\bf d} / d_+) (\Delta T)^2 = L ({\bf d})^2$.
\end{pf}

Lemma~\ref{lem:lem7} has an important consequence. Suppose that ${\mathcal A}$ is the
current active set for the Lasso, as in (\ref{eqn5.21}), and that ${\mathcal A}
\subseteq {\mathcal A}_{10}$. Then Lemma~\ref{lem:lem2} says that 
$L({\bf d})$ is~$\geq
A_{\mathcal A}$, and (\ref{eqn5.46}) gives 
\begin{equation}\label{eqn5.51}
\Delta S \geq -2 \widehat C \cdot \Delta T + A_{\mathcal A}^2 \cdot
(\Delta T)^2,
\end{equation}
with equality if {\bf d} is chosen to give the equiangular vector
${\bf u}_{\mathcal A}$, $d_{\mathcal A} = S_{\mathcal A} w_{\mathcal A}$, $d_{{\mathcal
    A}^c} = 0$. The Lasso operates to minimize $S(T)$ so we want
$\Delta S$ to be as negative as possible.  Lemma~\ref{lem:lem7} says that if the
support of {\bf d} is not confined to ${\mathcal A}_{10}$, 
then $\dot S(0)$
exceeds the optimum value $-2 \widehat C$; if it is confined, then
$\dot S(0) = -2 \widehat C$ but $\ddot S(0)$ exceeds the minimum value
$2 A_{\mathcal A}$ unless $d_{\mathcal A}$ is proportional to $S_{\mathcal A}
w_{\mathcal A}$ as in (\ref{eqn5.21}).

Suppose that $\widehat \bbeta$, a Lasso solution, exactly equals a
$\widehat \bbeta$ obtained from the Lasso-modified LARS algorithm,
henceforth called LARS--Lasso, as at $t = 1000$ in Figures 1 and 3.
We know from Lemma~\ref{lem:lem4} that subsequent Lasso estimates will follow a
linear track determined by some subset ${\mathcal A}$, $\bmu (\gamma) =
\widehat \bmu + \gamma {\bf u}_{\mathcal A}$, and so will the LARS--Lasso
estimates, but to verify Theorem 1 we need to show that ``${\mathcal A}$''
is the same set in both cases.

Lemmas 4--7 put four constraints on the Lasso 
choice of ${\mathcal A}$.
Define ${\mathcal A}_1$, ${\mathcal A}_0$ 
and ${\mathcal A}_{10}$ as at (\ref{eqn5.44}).

\begin{con} 
${\mathcal A}_1 \subseteq {\mathcal A}$. This
follows from Lemma~\ref{lem:lem4} since for sufficiently small $\gamma$ the
subsequent Lasso coefficients (\ref{eqn5.21}),
\begin{equation}
\widehat \beta_{\mathcal A} (\gamma) = \widehat \beta_{\mathcal A} + \gamma
S_{\mathcal A} w_{\mathcal A},
\end{equation}
will have $\widehat \beta_j (\gamma) \ne  0,  j \in {\mathcal
 A}_1$.
\end{con}

\begin{con} 
${\mathcal A} \subseteq {\mathcal A}_{10}$. 
Lemma~\ref{lem:lem7}, (\ref{eqn5.45}) shows that the Lasso 
choice~$\widehat {\bf d}$ in
$\bbeta(\gamma) = \widehat \bbeta + \gamma \widehat {\bf d}$ must have
its nonzero support in ${\mathcal A}_{10}$, or equivalently that
$\widehat \bmu (\gamma) = \widehat \bmu + \gamma {\bf u}_{\mathcal A}$
must have ${\bf u}_{\mathcal A} \in {\mathcal L}(X_{{\mathcal 
A}_{10}})$.  (It is
possible that ${\bf u}_{\mathcal A}$ happens to equal ${\bf u}_{\mathcal B}$
for some ${\mathcal B} \supset {\mathcal A}_{10}$, but that does not affect
the argument below.)
\end{con}

\begin{con} 
$w_{\mathcal A} = A_{\mathcal A} {\mathcal G}_{\mathcal
  A}^{-1} {\bf 1}_{\mathcal A}$  cannot have $\sign (w_j) \ne  
\sign (\widehat c_j)$ for any coordinate $j \in {\mathcal A}_0$. If
it does, then $\sign (\widehat \beta_j (\gamma)) \ne 
\mbox{sign} (\widehat c_j (\gamma))$ for sufficiently small~$\gamma$,
violating Lemma~\ref{lem:lem5}.
\end{con}

\begin{con} 
Subject to Constraints 1--3, ${\mathcal
 A}$ must minimize $A_{\mathcal A}$. This follows from 
Lemma~\ref{lem:lem7} as in (\ref{eqn5.51}),
and the requirement that the Lasso curve~$S(T)$ declines at the
fastest possible rate. 
\end{con}

Theorem 1 follows by induction: beginning at $\widehat \bbeta_0 = 0$,
we follow the LARS--Lasso algorithm and show that at every succeeding
step it must continue to agree with the Lasso definition (\ref{eqn1.5}). 
First
of all, 
suppose that $\widehat \bbeta$, our hypothesized Lasso and
LARS--Lasso solution, has occurred strictly {\it within} a LARS--Lasso
step. Then ${\mathcal A}_0$ is empty so that Constraints 1 and 2 imply
that ${\mathcal A}$ cannot change its current value: the equivalence
between Lasso and LARS--Lasso must continue at least to the end of the
step.

The one-at-a-time assumption of Theorem 1 says that at a LARS--Lasso
breakpoint, ${\mathcal A}_0$ 
has exactly one member, say $j_0$, 
so ${\mathcal A}$ must equal ${\mathcal A}_1$ or ${\mathcal A}_{10}$. 
There are two cases:
if $j_0$ has just been {\it added} to the set $\{ \vert \widehat c_j
\vert = \widehat C\}$, 
then Lemma~\ref{lem:lem1} says that $\mbox{sign}(w_{j_0}) =
\sign (\widehat c_{j_0})$, so that Constraint 3 
is not
violated; the other three constraints and Lemma~\ref{lem:lem2} imply that the Lasso
choice ${\mathcal A} = {\mathcal A}_{10}$ agrees with the LARS--Lasso
algorithm. The other case has $j_0$ {\it deleted} 
from the active set
as in~(\ref{eqn3.6}).  Now the choice ${\mathcal A} = {\mathcal A}_{10}$ is ruled out
by Constraint 3: it would keep $w_{\mathcal A}$ the same as in the
previous LARS--Lasso step, and we know that that was stopped in (\ref{eqn3.6})
to prevent a sign contradiction at coordinate $j_0$. In other words,
${\mathcal A} = {\mathcal A}_1$, in accordance with the Lasso modification of
LARS. This completes the proof of Theorem~1.

A LARS--Lasso algorithm is available even if the one-at-a-time
condition does not hold, but at the expense of additional computation.
Suppose, for example, 
{\it two} new members $j_1$ and $j_2$ are added to
the set $\{ \vert \widehat c_j \vert = \widehat C\}$, so ${\mathcal A}_0 =
\{j_1, j_2\}$. It is possible but not certain that ${\mathcal A}_{10}$
does not violate Constraint 3, in which case ${\mathcal A} = {\mathcal
  A}_{10}$. However, if it does violate Constraint 3, then both possibilities
${\mathcal A} = {\mathcal A}_1 \cup \{ j_1 \}$ and ${\mathcal A} = {\mathcal A}_1 \cup
\{ j_2\}$ must be examined to see which one gives the smaller value of
$A_{\mathcal A}$. Since one-at-a-time computations, perhaps with some
added~{\bf y} jitter, apply to all practical situations, the LARS
algorithm described in Section~\ref{sec:sec7} is not equipped to handle
many-at-a-time problems.

\section{Stagewise properties.}\label{sec:sec6} 
The main goal of this section
is to verify Theorem~2. Doing so also gives us a chance to make a more
detailed comparison of the LARS and Stagewise procedures. Assume that
$\widehat \bbeta$ is a Stagewise estimate of the regression
coefficients, for example, 
as indicated at $\sum \vert \widehat \beta_j
\vert = 2000$ in the right panel of Figure~\ref{fig1}, with prediction vector
$\widehat \bmu = X \widehat \bbeta$, current correlations $\widehat
{\bf c} = X^\prime ({\bf y} - \widehat \bmu)$, $\widehat C = \max\{
\vert \widehat c_j \vert \}$ 
and maximal set ${\mathcal A} = \{ j \dvtx \vert
\widehat c_j \vert = \widehat C \}$. We must show that successive
Stagewise estimates of $\bbeta$ develop according to the modified LARS
algorithm of Theorem 2, henceforth called LARS--Stagewise. For
convenience we can assume, by redefinition of ${\bf x}_j$ as $-{\bf
  x_j}$, if necessary, that the signs $s_j = \sign
  (\widehat c_j)$ are all non-negative.

As in (\ref{eqn3.8})--(\ref{eqn3.10}) we suppose that the 
Stagewise procedure (\ref{eqn1.7}) has
taken $N$~additional $\varepsilon$-steps forward from $\widehat \bmu = X
\widehat \bbeta$, giving new prediction vector $\widehat \bmu (N)$. 

\begin{lemma}
  \label{lem:lem8}
For sufficiently small $\varepsilon$, only $j
\in {\mathcal A}$ can have $P_j = N_j / N > 0$. 
\end{lemma}

\begin{pf}
Letting $N \varepsilon \equiv \gamma$, $\Vert
\widehat \bmu (N) - \widehat \bmu \Vert \leq \gamma$ so that $\widehat
{\bf c} (N) = X^\prime ({\bf y} - \widehat \bmu (N))$ satisfies
\begin{equation}
\vert \widehat c_j (N) - \widehat c_j \vert = 
\bigl\vert{\bf x}_j^\prime \bigl(
\widehat \bmu (N) - \widehat \bmu\bigr) \bigr\vert \leq \Vert {\bf x}_j \Vert
\cdot \Vert \widehat \bmu (N) - \widehat \bmu \Vert \leq \gamma .
\end{equation}
For $\gamma < \frac{1}{2} [\widehat C - \max_{{\mathcal A}^c} \{
\widehat c_j \}]$, $j$ in ${\mathcal A}^c$ cannot have maximal current
correlation and can never be involved in the $N$ steps.
\end{pf}

Lemma~\ref{lem:lem8} says that we can write the developing Stagewise prediction
vector as
\begin{equation}\label{eqn6.2}
\widehat \bmu (\gamma) = \widehat \bmu + \gamma {\bf v}, \qquad
\mbox{where }  {\bf v} = X_{\mathcal A} P_{\mathcal A},
\end{equation}
$P_{\mathcal A}$ a vector of length $\vert {\mathcal A} \vert$, with
components $N_j / N$ for $j \in {\mathcal A}$. The nature of the Stagewise
procedure puts three constraints on ${\bf v}$, the most obvious of
which is the following.

\begin{conI*} 
The vector ${\bf v} \in {\mathcal S}_{\mathcal A}^+$, the
nonnegative simplex 
\begin{equation}
{\mathcal S}_{\mathcal  A}^+ = \Biggl\{{\bf v} \dvtx  
{\bf v} = \sum_{j \in {\mathcal
 A}} {\bf x}_j P_j,  P_j \geq 0,  \sum_{j \in {\mathcal A}} P_j = 1
\Biggr\}. 
\end{equation}
Equivalently, $\gamma {\bf v} \in {\mathcal C}_{\mathcal A}$, the convex
cone (\ref{eqn3.12}).
\end{conI*}

The Stagewise procedure, unlike LARS, is not required to use all of
the maximal set ${\mathcal A}$ as the active set, and can instead restrict
the nonzero coordinates $P_j$ to a subset ${\mathcal B} \subseteq {\mathcal
  A}$. Then ${\bf v} \in {\mathcal L}(X_{\mathcal B})$, the linear space
spanned by the columns of $X_{\mathcal B}$, but not all such vectors ${\bf
  v}$ are allowable Stagewise forward directions.

\begin{conII*} 
The vector ${\bf v}$ must be proportional to the
equiangular vector ${\bf u}_{\mathcal B}$, (\ref{eqn2.6}), that is, 
${\bf v} = {\bf
 v}_{\mathcal B}$, (\ref{eqn5.12}),
\begin{equation}\label{eqn6.4}
{\bf v}_{\mathcal B} = A_{\mathcal B}^2 X_{\mathcal B} {\mathcal G}_{\mathcal B}^{-1}
{\bf 1}_{\mathcal B} = A_{\mathcal B} {\bf u}_{\mathcal B}.
\end{equation}
\end{conII*}

Constraint II amounts to requiring that the current correlations in
${\mathcal B}$ decline at an equal rate: since
\begin{equation}\label{eqn6.5}
\widehat c_j (\gamma) = {\bf x}_j^\prime ({\bf y} - \widehat \bmu -
\gamma {\bf v}) = \widehat c_j - \gamma {\bf x}_j^\prime {\bf v},
\end{equation}
we need $X_{\mathcal B}^\prime {\bf v} = \lam {\bf 1}_{\mathcal B}$ for
some $\lam > 0$, implying ${\bf v} = \lam {\mathcal G}_{\mathcal B}^{-1}
{\bf 1}_{\mathcal B}$; choosing $\lam = A_{\mathcal B}^2$ satisfies Constraint
II. Violating Constraint II makes the current correlations $\widehat
c_j(\gamma)$ unequal so that the Stagewise algorithm as defined at
(\ref{eqn1.7}) could not proceed in direction ${\bf v}$.

Equation (\ref{eqn6.4}) gives $X_{\mathcal B}^\prime {\bf v}_{\mathcal B} = A_{\mathcal
  B}^2 {\bf 1}_{\mathcal B}$, or
\begin{equation}\label{eqn6.6}
{\bf x}_j^\prime {\bf v}_{\mathcal B} = A_{\mathcal B}^2 \qquad 
\mbox{for }
j \in {\mathcal B}. 
\end{equation}

\begin{conIII*} 
The vector ${\bf v} = {\bf v}_{\mathcal
  B}$ must satisfy
\begin{equation}\label{eqn6.7}
{\bf x}_j^\prime {\bf v}_{\mathcal B} \geq A_{\mathcal B}^2 \qquad 
\mbox{for }
j \in {\mathcal A} - {\mathcal B}.
\end{equation}
\end{conIII*}

Constraint III follows from (\ref{eqn6.5}). It says that the current
correlations for members of ${\mathcal A} = \{ j \dvtx \vert \widehat c_j
\vert = \widehat C\}$ {\it not} in ${\mathcal B}$ must decline at least as
quickly as those in~${\mathcal B}$. If this were not true, then ${\bf
  v}_{\mathcal B}$ would not be an allowable direction for Stagewise
development since variables in ${\mathcal A} - {\mathcal B}$ would immediately
reenter (\ref{eqn1.7}).

To obtain strict inequality in (\ref{eqn6.7}), let ${\mathcal B}_0 
\subset {\mathcal A}
- {\mathcal B}$ be the set of indices for which ${\bf x}_j^\prime {\bf
  v}_{\mathcal B} = A_{\mathcal B}^2$. It is easy to show that ${\bf v}_{{\mathcal
    B} \cup {\mathcal B}_o} = {\bf v}_{\mathcal B}$. In other words, 
if we 
take~${\mathcal B}$ to be the {\it largest} set having a given ${\bf v}_{\mathcal
  B}$ proportional to its equiangular vector, then ${\bf x}_j^\prime
{\bf v}_{\mathcal B} > A_{\mathcal B}^2$ for 
$j \in {\mathcal A} - {\mathcal B}$.

Writing $\widehat \bmu (\gamma) = \widehat \bmu + \gamma {\bf v}$ as
in (\ref{eqn6.2}) presupposes that the Stagewise solutions follow a piecewise
linear track. However, the presupposition can be reduced to one of
piecewise differentiability by taking $\gamma$ infinitesimally small.
We can always express the family of Stagewise solutions as $\widehat
\bbeta(z)$, where the real-valued parameter~$Z$ plays the role of $T$
for the Lasso, increasing from 0 to some maximum value as $\widehat
\bbeta(z)$ goes from ${\bf 0}$ to the full OLS estimate. [The
choice $Z = T$ used in Figure~\ref{fig1} may not necessarily yield a one-to-one
mapping; $Z = S({\bf 0}) - S(\widehat \bbeta)$, the reduction in
residual squared error, always does.] We suppose that the Stagewise
estimate $\widehat {\bbeta}(z)$ is everywhere right differentiable
with respect to $z$. Then the right derivative
\begin{equation}
\widehat {\bf v} = d \widehat \bbeta (z) / d z
\end{equation}
must obey the three constraints.

The definition of the idealized Stagewise procedure in Section~\ref{sec:sec3.2}, in
which \mbox{$\varepsilon \rightarrow 0$} in rule (\ref{eqn1.7}), is somewhat vague but
the three constraints apply to any reasonable interpretation. It turns
out that the LARS--Stagewise algorithm satisfies the constraints and is
unique in doing so. This is the meaning of Theorem 2. [Of course the
LARS--Stagewise algorithm is also supported by direct numerical
comparisons with (\ref{eqn1.7}), as in Figure~\ref{fig1}'s right panel.]

If ${\bf u}_{\mathcal A} \in {\mathcal C}_{\mathcal A}$, then ${\bf v} = {\bf
  v}_{\mathcal A}$ obviously satisfies the three constraints. The
interesting situation for Theorem 2 is ${\bf u}_{\mathcal A} \notin
{\mathcal C}_{\mathcal A}$, which we now assume to be the case. Any subset
${\mathcal B} \subset {\mathcal A}$ determines a face of the convex cone of
dimension~$\vert {\mathcal B} \vert$, the face having $P_j > 0$ in (\ref{eqn3.12})
for $j \in {\mathcal B}$ and $P_j = 0$ for $j \in {\mathcal A} - {\mathcal B}$.
The orthogonal projection of ${\bf u}_{\mathcal A}$ into the linear
subspace ${\mathcal L}(X_{\mathcal B})$, say $\mbox{Proj}_{{\mathcal B}}({\bf
  u}_{\mathcal A})$, is proportional to ${\mathcal B}$'s equiangular vector
${\bf u}_{\mathcal B}$: using (\ref{eqn2.7}),
\begin{equation}\label{eqn6.9}
\Proj_{\mathcal B} ({\bf u}_{\mathcal A}) = X_{\mathcal B} {\mathcal G}_{\mathcal
  B}^{-1} X_{\mathcal B}^\prime {\bf u}_{\mathcal A} = X_{\mathcal B} {\mathcal
  G}_{\mathcal B}^{-1} A_{\mathcal A} {\bf 1}_{\mathcal B} = (A_{\mathcal A} / A_{\mathcal B})
  \cdot {\bf u}_{\mathcal B},
\end{equation}
or equivalently
\begin{equation}\label{eqn6.10}
\Proj_{\mathcal B}({\bf v}_{\mathcal A}) = (A_{\mathcal A} / A_{\mathcal B})^2
{\bf v}_{\mathcal B}.
\end{equation}

\begin{figure}[h]

\includegraphics[scale=1.1]{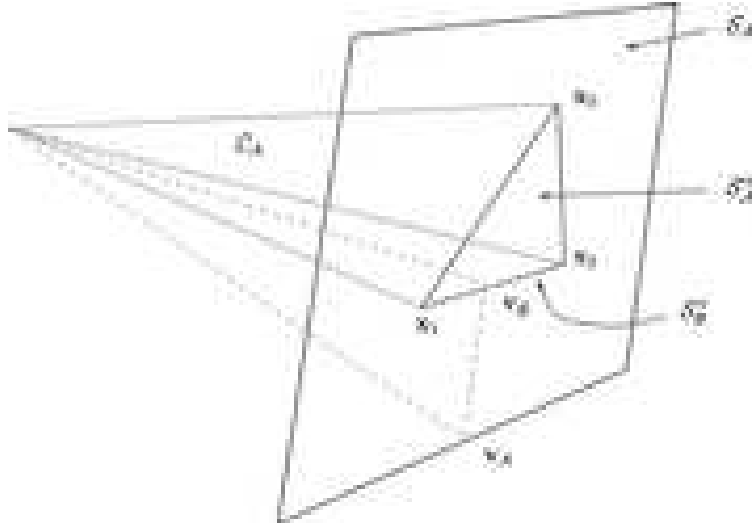}

\caption{The geometry of the
    LARS--Stagewise modification.}
\label{fig8}
\end{figure}

The nearest point to ${\bf u}_{\mathcal A}$ in ${\mathcal C}_{\mathcal A}$, say
$\widehat {\bf u}_{\mathcal A}$, 
is of the form $\Sig_{\mathcal A} {\bf x}_j
\widehat P_j$ with $\widehat P_j \geq 0$. Therefore $\widehat {\bf
  u}_{\mathcal A}$ exists strictly within face $\widehat {\mathcal B}$, where
$\widehat {\mathcal B} = \{j \dvtx \widehat P_j > 0\}$, and must equal
$\mbox{Proj}_{\hat {\mathcal B}} ({\bf u}_{\mathcal A})$. According to (\ref{eqn6.9}),
$\widehat {\bf u}_{\mathcal A}$ is proportional to $\widehat {\mathcal B}$'s
equiangular vector ${\bf u}_{\hat {\mathcal B}}$, and also to ${\bf
  v}_{\hat {\mathcal B}} = A_{\mathcal B} {\bf u}_{\mathcal B}$. In other words
${\bf v}_{\hat {\mathcal B}}$ satisfies Constraint II, and it obviously
also satisfies Constraint I. Figure~\ref{fig8} schematically illustrates the
geometry.

\begin{lemma}
\label{lem:lem9}
The vector ${\bf v}_{\hat {\mathcal B}}$ satisfies
Constraints \textup{I--III}, and conversely if ${\bf v}$ satisfies the three
constraints, then ${\bf v} = {\bf v}_{\hat {\mathcal B}}$.
\end{lemma}

\begin{pf}
 Let $\Cos \equiv A_{\mathcal A} / A_{\mathcal B}$
and $\Sin = [1 - \mbox{Cos}^2]^{1/2}$, the latter being
greater than zero by Lemma~\ref{lem:lem2}. For any face ${\mathcal B} \subset {\mathcal
  A}$, (\ref{eqn6.9}) implies 
\begin{equation}\label{eqn6.11}
{\bf u}_{\mathcal A} = \Cos \cdot {\bf u}_{\mathcal B} + \Sin
\cdot {\bf z}_{\mathcal B},
\end{equation}
where ${\bf z}_{\mathcal B}$ is a unit vector orthogonal to ${\mathcal
  L}(X_{\mathcal B})$, pointing away from ${\mathcal C}_{\mathcal A}$. By an
$n$-dimensional coordinate rotation we can make 
${\mathcal L}(X_{\mathcal B})
= {\mathcal L}({\bf c}_1, {\bf c}_2, \ldots, {\bf c}_J)$, $J = \vert {\mathcal
  B} \vert$, the space of $n$-vectors with last $n-J$ coordinates
zero, and also
\begin{equation}\label{eqn6.12}
{\bf u}_{\mathcal B} = (1, {\bf 0}, 0, {\bf 0}), \qquad 
{\bf u}_{\mathcal A} =
(\Cos, {\bf 0}, \Sin, {\bf 0}),
\end{equation}
the first ${\bf 0}$ having length $J-1$, the second {\bf 0}
length $n - J - 1$. Then we can write
\begin{equation}
{\bf x}_j = \bigl(A_{\mathcal B}, {\bf x}_{j_2}, 0, {\bf 0}\bigr) \qquad 
\mbox{for }
j \in {\mathcal B},
\end{equation}
the first coordinate $A_{\mathcal B}$ being required since ${\bf
  x}_j^\prime {\bf u}_{\mathcal B} = A_{\mathcal B}$, (\ref{eqn2.7}). Notice that ${\bf
  x}_j^\prime {\bf u}_{\mathcal A} = \Cos 
\cdot A_{\mathcal B} = A_{\mathcal
  A}$, as also required by (\ref{eqn2.7}).

For $\ell \in {\mathcal A} - {\mathcal B}$ denote ${\bf x}_{\ell}$ as
\begin{equation}
{\bf x}_{\ell} = \bigl(x_{\ell_1}, {\bf x}_{\ell_2}, x_{\ell_3}, {\bf
  x}_{\ell_4}\bigr),
\end{equation}
so (\ref{eqn2.7}) yields
\begin{equation}\label{eqn6.15}
A_{\mathcal A} = {\bf x}_{\ell}^\prime {\bf u}_{\mathcal A} = \Cos \cdot
x_{\ell_1} + \Sin \cdot x_{\ell_3} .
\end{equation}
Now assume ${\mathcal B} = \widehat {\mathcal B}$. In this case a
separating hyperplane ${\mathcal H}$ orthogonal 
to ${\bf z}_{\hat {\mathcal
    B}}$~in~(\ref{eqn6.11}) passes between the convex cone ${\mathcal C}_{\mathcal A}$
and ${\bf u}_{\mathcal A}$, through $\widehat {\bf u}_{\mathcal A} =
\Cos \cdot {\bf u}_{\hat {\mathcal B}}$, implying $x_{\ell_3} \leq
0$ [i.e., ${\bf x}_{\ell}$ and ${\bf u}_{\mathcal A}$ are on opposite
sides of ${\mathcal H}$, $x_{\ell_3}$ being negative since the
corresponding coordinate of ${\bf u}_{\mathcal A}$, ``Sin'' in (\ref{eqn6.12}), is
positive]. Equation~(\ref{eqn6.15}) gives $\Cos \cdot x_{\ell_1} \geq
A_{\mathcal A} = \Cos \cdot A_{\hat {\mathcal B}}$ or
\begin{equation}
{\bf x}_{\ell}^\prime {\bf v}_{\hat {\mathcal B}} = x_{\ell}^\prime
(A_{\hat {\mathcal B}} {\bf u}_{\hat {\mathcal B}}) = A_{\hat {\mathcal B}}
x_{\ell_1} \geq A_{\hat {\mathcal B}}^2,
\end{equation}
verifying that Constraint III is satisfied.

Conversely suppose that ${\bf v}$ satisfies Constraints I--III so that
${\bf v} \in {\mathcal S}_{\mathcal A}^+$ and ${\bf v} = {\bf v}_{\mathcal B}$ for
the nonzero coefficient set ${\mathcal B}$: ${\bf v}_{\mathcal B} =
\Sig_{\mathcal B} {\bf x}_j P_j, P_j > 0$. Let ${\mathcal H}$ be the
hyperplane passing through $\mbox{Cos} \cdot {\bf u}_{\mathcal B}$
orthogonally to ${\bf z}_{\mathcal B}$, (\ref{eqn6.9}), (\ref{eqn6.11}). 
If ${\bf v}_{\mathcal
  B}  \ne  {\bf v}_{\hat {\mathcal B}}$, 
then at least one of the
vectors ${\bf x}_{\ell}$, $\ell \in {\mathcal A} - {\mathcal B}$, must lie on
the same side of ${\mathcal H}$ as ${\bf u}_{\mathcal A}$, so that~$x_{\ell_3}
> 0$ (or else ${\mathcal H}$ would be a separating hyperplane between
${\bf u}_{\mathcal A}$ and ${\mathcal C}_{\mathcal A}$, and ${\bf v}_{\mathcal B}$
would be proportional to $\widehat {\bf u}_{\mathcal A}$, the nearest
point to ${\bf u}_{\mathcal A}$ in ${\mathcal C}_{\mathcal A}$, implying ${\bf
  v}_{\mathcal B} = {\bf v}_{\hat {\mathcal B}}$). Now (\ref{eqn6.15}) gives
$\mbox{Cos} \cdot x_{\ell_1} < A_{\mathcal A} = \mbox{Cos} \cdot A_{\mathcal
B}$, or
\begin{equation}
{\bf x}_{\ell}^\prime {\bf v}_{\mathcal B} = {\bf x}_{\ell}^\prime
(A_{\mathcal B} {\bf u}_{\mathcal B}) = A_{\mathcal B} x_{\ell_1} < A_{\mathcal B}^2 .
\end{equation}
This violates Constraint III, showing that ${\bf v}$ must equal
${\bf v}_{\hat {\mathcal B}}$. 
\end{pf}

Notice that the direction of advance $\widehat {\bf v} = {\bf v}_{\hat
  {\mathcal B}}$ of the idealized Stagewise procedure is a function only
of the current maximal set $\widehat {\mathcal A} = \{ j \dvtx \vert \widehat
c_j \vert = \widehat C\}$, say $\widehat {\bf v} = \phi (\widehat
{\mathcal A}\,)$. In the language of (\ref{eqn6.7}),
\begin{equation}\label{eqn6.18}
\frac{d \widehat \bbeta(z)}{dz} = \phi (\widehat {\mathcal A}\,).
\end{equation}

The LARS--Stagewise algorithm of Theorem 2 produces an evolving family
of estimates $\widehat \bbeta$ that everywhere satisfies (\ref{eqn6.18}). This
is true at every LARS--Stagewise breakpoint by the definition of the
Stagewise modification. It is also true between breakpoints. Let
$\widehat {\mathcal A}$ be the maximal set at the breakpoint, giving
$\widehat {\bf v} = {\bf v}_{\hat {\mathcal B}} = \phi (\widehat {\mathcal
  A})$. In the succeeding LARS--Stagewise interval $\widehat \bmu
(\gamma) = \widehat \bmu + \gamma {\bf v}_{\hat {\mathcal B}}$, the
maximal set is immediately reduced to $\widehat {\mathcal B}$, according
to properties (\ref{eqn6.6}), (\ref{eqn6.7}) of ${\bf v}_{\hat {\mathcal 
B}}$, 
at which it
stays during the entire interval. However, $\phi (\widehat {\mathcal 
B}\,) =
\phi (\widehat {\mathcal A}\,) = {\bf v}_{\hat {\mathcal B}}$ 
since ${\bf
  v}_{\hat {\mathcal B}} \in {\mathcal C}_{\hat {\mathcal B}}$, so the
LARS--Stagewise procedure, which continues in the direction $\widehat
{\bf v}$ until a new member is added to the active set, continues to
obey the idealized Stagewise equation (\ref{eqn6.18}).

All of this shows that the LARS--Stagewise algorithm produces a
legitimate version of the idealized Stagewise track. The converse of
Lemma~\ref{lem:lem9} says that there are no other versions, verifying Theorem 2.

The Stagewise procedure has its potential generality as an advantage
over LARS and Lasso: it is easy to define forward Stagewise methods
for a wide variety of nonlinear fitting problems, as in 
\citeauthor{ElemStatLearn} [(\citeyear{ElemStatLearn}), Chapter~10, 
which begins with a Stagewise analysis 
of
``boosting'']. Comparisons with LARS and Lasso within the linear model
framework, as at the end of Section~\ref{eqn3.2}, help us better understand
Stagewise methodology. This section's results permit further
comparisons.

Consider proceeding forward from $\widehat \bmu$ along unit vector
${\bf u}$, $\widehat \bmu (\gamma) = \widehat \bmu + \gamma {\bf u}$,
two interesting choices being the LARS direction ${\bf u}_{\hat {\mathcal
    A}}$ and the Stagewise direction $\widehat \bmu_{\hat {\mathcal B}}$.
For ${\bf u} \in {\mathcal L}(X_{\hat {\mathcal A}})$, the rate of change of
$S(\gamma) = \Vert {\bf y} - \widehat \bmu (\gamma) \Vert^2$ is
\begin{equation}\label{eqn6.19}
-\frac{\partial S(\gamma)}{\partial \gamma} \bigg|_0 = 2 \widehat C
  \cdot \frac{{\bf u}_{\mathcal A}^\prime \cdot {\bf u}}{A_{\hat {\mathcal
  A}}},
\end{equation}
(\ref{eqn6.19}) following quickly from (\ref{eqn5.18}). This shows that the LARS
direction ${\bf u}_{\hat {\mathcal A}}$ maximizes the instantaneous
decrease in $S$. The ratio
\begin{equation}\label{eqn6.20}
\frac{\partial S_{\rm Stage} (\gamma)}{\partial \gamma} \bigg|_0
\,\Big/\, \frac{\partial S_{\rm LARS} (\gamma)}{\partial \gamma}
\bigg|_0 = \frac{A_{\hat {\mathcal A}}}{A_{\hat {\mathcal B}}},
\end{equation}
equaling the quantity ``Cos'' in (\ref{eqn6.15}).

The comparison goes the other way for the maximum absolute correlation
$\widehat C(\gamma)$. Proceeding as in (\ref{eqn2.15}), 
\begin{equation}\label{eqn6.21}
- \frac{\partial \widehat C(\gamma)}{\partial \gamma} \bigg|_0 =
  \min_{\hat {\mathcal A}} \{ \vert x_j^\prime {\bf u} \vert \}.
\end{equation}
The argument for Lemma~\ref{lem:lem9}, using Constraints II and III, shows
that ${\bf u}_{\hat {\mathcal B}}$ maximizes (\ref{eqn6.21}) 
at $A_{\hat {\mathcal
    B}}$, and that
\begin{equation}\label{eqn6.22}
\frac{\partial \widehat C_{\rm LARS}(\gamma)}{\partial \gamma}
\bigg|_0  \,\Big/\, \frac{\partial \widehat C_{\rm
    Stage}(\gamma)}{\partial \gamma} \bigg|_0 
= \frac{A_{\hat {\mathcal
     A}}}{A_{\hat {\mathcal B}}}.
\end{equation}

The original motivation for the Stagewise procedure was to minimize
residual squared error within a framework of parsimonious forward
search. However, (\ref{eqn6.20}) shows that Stagewise is less greedy than LARS
in this regard, it being more accurate to describe Stagewise as
striving to minimize the maximum absolute residual correlation.

\section{Computations.}\label{sec:sec7}
The entire sequence of steps in the LARS algorithm with~$m<n$
variables requires $O(m^3 +nm^2)$ computations---the cost of a least
squares fit on $m$ variables. 

In detail, at the $k$th of $m$ steps, we compute $m-k$ inner products $c_{jk}$
of the nonactive ${\bf x}_j$ with the current residuals to identify
the next active variable, and then invert the $k\times k$ matrix $\G_k=X_k'X_k$
to find the next LARS direction.  We do this by updating the Cholesky
factorization $R_{k-1}$ of $\G_{k-1}$ found at the previous step
[\citet{GVL83}]. At the final step $m$, we have computed the Cholesky
$R=R_m$ for the full cross-product matrix, which is the dominant
calculation for a least squares fit. Hence the LARS sequence can be
seen as a Cholesky factorization with a guided ordering of the
variables. 

The computations can be reduced further by recognizing that the inner
products above can be updated at each iteration using the
cross-product matrix $X'X$ and the current directions. For $m\gg n$,
this strategy is counterproductive and is not used.

For the {\em lasso} modification, the computations are similar, except
that occasionally one has to drop a variable, and hence {\em downdate}
$R_k$ [costing at most $O(m^2)$ operations per downdate]. 
For the {\em stagewise}
modification of LARS, we need to check at each iteration that the
components of $w$ are all positive.  If not, one or more variables are
dropped [using the {\em inner loop} of the NNLS algorithm described in
\citet{LH74}], again requiring downdating of $R_k$.  With 
many
correlated variables, the stagewise version can take many more steps
than LARS because of frequent dropping and adding of variables,
increasing the computations by a factor up to 5 or more in extreme
cases.

The LARS algorithm (in any of the three states above) works
gracefully for the case where there are many more variables than
observations: $m\gg n$. In this case
LARS terminates at the saturated least squares fit after $n-1$ variables
have entered the active set [at a cost of $O(n^3)$ operations]. (This
number is $n-1$ rather than $n$, because the columns of $X$ have been
mean centered, and hence it has row-rank $n-1$.) 
We make a few more remarks
about the $m\gg n$ case in 
the {\em lasso} state:
\begin{enumerate}
\item The LARS algorithm continues to provide Lasso solutions along
the way, and the final solution highlights the fact that a Lasso fit
can have no more than $n-1$ (mean centered) variables with nonzero
coefficients.

\item Although the model involves no more than $n-1$ variables at any
  time, the number of {\em different} variables ever to have entered
  the model during the entire sequence can be---and typically is---greater
  than $n-1$.

\item The model sequence, particularly near the saturated end, tends
  to be quite variable with respect to small changes in $\bf y$.

\item The estimation of $\sigma^2$ may have to depend on an auxiliary
  method such as nearest neighbors (since the final model is
  saturated). We have not investigated the accuracy of the simple
  approximation formula (\ref{eq:stein}) for the case $m>n$.
\end{enumerate}

Documented S-PLUS implementations of LARS and associated functions are
available from www-stat.stanford.edu/$\sim$hastie/Papers/; the
diabetes data also appears there.

\section{Boosting procedures.}\label{sec:sec8}
One motivation for studying the Forward Stagewise algorithm is its
usefulness in adaptive fitting for data mining. In particular, Forward
Stagewise
ideas are used in ``boosting,'' an
 important class of 
fitting methods for data mining introduced by \citet{FS95}.
These methods are one of the hottest topics in the area of  machine
learning, and one of the most effective prediction methods in current use.
Boosting can use any adaptive
fitting procedure as its ``base learner'' (model fitter):
 trees are a popular choice, as implemented in
CART [\citet{BFOS84}].

\citet{FHT00} and \citet{Fr99} 
studied boosting and proposed a number of
procedures,
the most relevant to this discussion being {\em least squares boosting}.
This procedure works by 
successive fitting of regression trees to the current  residuals.
Specifically we start with the residual $\br=\by$ and the fit $\hat
\by = 0$.
We  fit a tree in ${\bx_1, \bx_2, \ldots, \bx_m}$ to
 the response ${\bf y}$ giving a fitted tree~${\bf t}_1$ (an~$n$-vector of
fitted values).
Then we update $  \hat\by$ to $\hat\by+ \varepsilon\cdot  {\bf t}_1 $,
$ \br$ to $\by-\hat\by $ and continue for many iterations.
Here $\varepsilon$ is a small positive constant. Empirical studies 
show
that small values of $\varepsilon$ work better than $\varepsilon=1$: in fact, 
for prediction accuracy ``the smaller the better.'' The only drawback
in taking very small values of $\varepsilon$ is  computational slowness.

A major research question has been why boosting works so well, and
specifically why is $\varepsilon$-shrinkage so important?  To
understand boosted trees in the present context, we think of our
predictors not as our original variables ${\bx_1, \bx_2, \ldots,
  \bx_m}$, but instead as the set of all trees ${\bf t}_k$ that could
be fitted to our data.  There is a strong similarity between
least squares boosting and Forward Stagewise regression as defined
earlier. Fitting a tree to the current residual is a numerical way of
finding the ``predictor'' most correlated with the residual. Note,
however, that the greedy algorithms used in CART do not search among
all possible trees, but only a subset of them. In addition the set
of all trees, including a parametrization for the predicted values in
the terminal nodes, is infinite. Nevertheless one can define idealized versions
of least-squares boosting that look much like Forward Stagewise
regression.

\citet{ElemStatLearn}
 noted the the striking similarity between Forward Stagewise regression
and the Lasso, and conjectured that this may help explain the success
of the  Forward Stagewise process used in least squares boosting.
That is, in some sense least squares boosting may be carrying out a
Lasso fit on the infinite set of tree predictors. 
Note that direct computation of the Lasso via  the LARS procedure 
would not be feasible
in this setting because the number of trees is infinite and  one 
could not compute the optimal step length. However, Forward Stagewise regression is
feasible because it only need find the the most correlated predictor
among the infinite set, where it approximates by  numerical
search.

In this paper we have established the connection between the Lasso and Forward
Stagewise regression.
We are now thinking about how these results can help to understand and
improve boosting procedures. One such idea is a modified form of
Forward Stagewise: we find the best tree as usual, but rather than
taking a small step in only that tree, we take a small least squares
step in all trees  currently in our model. One can show that for
small step sizes this procedure approximates LARS; its advantage is
that it can be carried out on an infinite set of predictors such as
trees.

\begin{appendix}
\appsection*{Appendix}
\label{sec:appendix}

\section{Local linearity and Lemma~\ref{prop:divergence}.}\label{sec:local-line-prop}

\renewcommand{\theequation}{\Alph{section}.\arabic{equation}}

\begin{conv*}
We write ${\bf x}_l$ with subscript $l$ for members
of the active set~${\mathcal A}_k$. Thus ${\bf x}_l$ denotes the $l$th variable
to enter, being an abuse of notation for $s_l {\bf x}_{j(l)} = \mbox{sgn}
(\hat c_{j(l)}) {\bf x}_{j(l)}$.
Expressions ${\bf x}_l'( {\bf y} - \hat \bmu_{k-1}({\bf y})
) = \hat C_k ({\bf y})$ and ${\bf x}_l'{\bf u}_k = A_k$ clearly do not depend on which
${\bf x}_l \in {\mathcal A}_k$ we choose.

By writing $j \notin {\mathcal A}_k$, we intend that both ${\bf x}_j$ and $-{\bf x}_j$
are candidates for inclusion at the next step. One could think of
negative indices $-j$ corresponding to ``new'' variables ${\bf x}_{-j} =-
{\bf x}_j$. 

The active set ${\mathcal A}_k ({\bf y})$ depends on the data ${\bf y}$. When ${\mathcal
  A}_k ({\bf y})$ is the same for all ${\bf y}$~in a neighborhood of ${\bf y}_0$, we 
  say that ${\mathcal A}_k ({\bf y})$ is locally fixed 
[at ${\mathcal A}_k = {\mathcal
  A}_k ({\bf y}_0)$].

A function $g({\bf y})$ is locally Lipschitz at ${\bf y}$ if for all
sufficiently small vectors $\Delta {\bf y}$, 
\begin{equation}
  \label{eq:lipdef}
 \Vert \Delta g \Vert  = \Vert g({\bf y} + \Delta 
{\bf y}) - g({\bf y}) \Vert  \leq L \Vert \Delta {\bf y} \Vert.
\end{equation}
If the constant $L$ applies for all ${\bf y}$, we say that $g$ is uniformly
locally Lipschitz $(L)$, and the word ``locally'' may be dropped.
\end{conv*}

\begin{lemma} \label{lem-gk}
  For each $k$, $0 \leq k \leq m$, there is an open set $G_k$ of full
  measure on which ${\mathcal A}_k({\bf y})$ and ${\mathcal A}_{k+1}({\bf y})$ are locally
  fixed and differ by $1$, and $\hat \bmu_k ({\bf y})$ is locally linear.
  The sets $G_k$ are decreasing as $k$ increases.
\end{lemma}

\begin{pf}
The argument is by induction.
The induction hypothesis states 
that for each ${\bf y}_0 \in G_{k-1}$ there is a small ball $B({\bf y}_0)$ on
which (a) the active sets
${\mathcal A}_{k-1} ({\bf y})$ and ${\mathcal A}_k ({\bf y})$ are fixed and equal to $
{\mathcal A}_{k-1}$ and ${\mathcal A}_k$, respectively, 
(b) $|{\mathcal A}_k \setminus {\mathcal A}_{k-1}| = 1$ so that the same
single variable enters locally at stage $k-1$ and (c)
$\hat \bmu_{k-1} ({\bf y}) = M{\bf y}$ is linear. We construct a set
$G_k$ with the same property.

Fix a point ${\bf y}_0$ and the corresponding ball $B({\bf y}_0)
\subset G_{k-1}$, on which ${\bf y}  - \hat \bmu_{k-1} ({\bf y})  = {\bf y} - M {\bf y} = R {\bf y}$, say.
 For indices $j_1, j_2 \notin {\mathcal A}$, let 
$N(j_1,j_2)$ be the set of ${\bf y}$ for which there exists a $\gamma$ such
that 
\begin{equation}
  \label{eq:double}
  w'(R{\bf y} - \gamma {\bf u}_k) = {\bf x}_{j_1}'( R{\bf y} - \gamma {\bf u}_k ) 
                      = {\bf x}_{j_2}'( R{\bf y} - \gamma {\bf u}_k ). 
\end{equation}
Setting $\bdelta_1 = {\bf x}_l - {\bf x}_{j_1}$, the first equality may be written
$\bdelta_1'R{\bf y} = \gamma \bdelta_1'{\bf u}_k$ and so when 
$\bdelta_1' {\bf u}_k \neq 0$
determines
\[
\gamma = \bdelta_1'R{\bf y}  / \bdelta_1'{\bf u}_k =: \bfeta_1'{\bf y}.
\]
[If $\bdelta_1' {\bf u}_k = 0$, there are no qualifying ${\bf y}$, and $N(j_1,j_2)$
is empty.]
Now using the second equality and setting $\bdelta_2 = {\bf x}_l - {\bf x}_{j_2}$, 
we see that $N(j_1, j_2)$ is contained
in~the set of ${\bf y}$ for which 
\[
  \bdelta_2'R{\bf y} = \eta_1'{\bf y}  \ \bdelta_2'{\bf u}_k.
\]
In other words, setting $\bfeta_2 = R'\bdelta_2 - (\bdelta_2'{\bf u}_k) \bfeta_1$,
we have
\[
  N(j_1,j_2) \subset \{{\bf y}\dvtx \bfeta_2'{\bf y} = 0 \}.
\]
If we define
\[
  N({\bf y}_0) = \bigcup \{ N(j_1,j_2)\dvtx j_1, j_2 \notin 
{\mathcal A}, j_1 \neq
  j_2 \},
\]
it is evident that $N({\bf y}_0)$ is a finite union of hyperplanes and hence
closed. For ${\bf y} \in B({\bf y}_0) \setminus N({\bf y}_0)$, a unique new variable
joins the active set at step $k$. Near each such ${\bf y}$ the ``joining''
variable is locally the same and $\gamma_k({\bf y}) {\bf u}_k$ is locally linear.

We then define $G_k \subset G_{k-1}$ 
as the union of such sets $B({\bf y}) \setminus N({\bf y})$
over ${\bf y} \in G_{k-1}$. Thus $G_k$ is open and, on $G_k$, ${\mathcal
  A}_{k+1}({\bf y})$ is locally constant and 
$\hat \bmu_k({\bf y})$ is locally
linear. Thus properties (a)--(c) hold for $G_k$.

The same argument works for the initial case $k=0$: since $\hat \bmu_0
= 0$, there is no circularity.

Finally, since the intersection of $G_k$ with any compact set is
covered by a finite number of $B(y_i) \setminus N(y_i)$, it is clear
that $G_k$ has full measure.
\end{pf}

\begin{lemma} \label{lem:singlepoint}
  Suppose that, for ${\bf y}$ near ${\bf y}_0$, $\hat \bmu_{k-1} ({\bf y})$ is
  continuous \textup{(}resp. linear\textup{)} 
and that ${\mathcal A}_k ({\bf y}) = {\mathcal A}_k$.
  Suppose also that, at ${\bf y}_0$, ${\mathcal A}_{k+1}({\bf y}_0) = {\mathcal A} \cup \{k+1\}$.

 Then  for ${\bf y}$ near ${\bf y}_0$,  ${\mathcal A}_{k+1}({\bf y}) = {\mathcal A}_k \cup
 \{ k+1 \}$ and $\hat \gamma_k({\bf y})$ and hence $\hat \bmu_k ({\bf y})$
are continuous \textup{(}resp. linear\textup{)} and
 uniformly Lipschitz.
\end{lemma}

\begin{pf}
  Consider first the situation at ${\bf y}_0$, with $\widehat C_k$ and
  $\widehat c_{kj}$ defined in (\ref{eqn2.18})~and~(\ref{eqn2.17}), respectively.
Since $k+1 \notin {\mathcal A}_k$, we have $|\hat C_k({\bf y}_0)| > \hat
c_{k,k+1} ({\bf y}_0)$, and $\hat \gamma_k ({\bf y}_0) > 0 $ satisfies
\begin{equation}
      \label{eq:singlept}
\qquad \hspace*{3pt}\hat C_k ({\bf y}_0) - \hat \gamma_k ({\bf y}_0) A_k \ 
    \left\{  \matrix{ = \cr  > } \right\} \ 
    \hat c_{k,j} ({\bf y}_0) - \hat \gamma_k ({\bf y}_0) a_{k,j} 
    \qquad \mbox{as } 
    \left\{  \matrix{ j = k+1 \cr j > k+1}.\right.
\end{equation}
In particular, it must be that $A_k \neq a_{k,k+1}$, 
and hence
\[
  \hat \gamma_k({\bf y}_0) = \frac{\hat C_k({\bf y}_0) - 
\hat c_{k,k+1}({\bf y}_0)}{A_k
  - a_{k,k+1}} > 0.
\]

Call an index $j$ admissible if $j \notin {\mathcal A}_k$ and $a_{k,j}
\neq A_k$. For ${\bf y}$ near ${\bf y}_0$, this property is independent of ${\bf y}$.
For admissible $j$, define
\[
  R_{k,j}({\bf y}) = \frac{\hat C_k({\bf y}) - \hat c_{k,j} ({\bf y})}{ A_k - a_{k,j}},
\]
which is continuous (resp. linear) near ${\bf y}_0$ from the assumption on
$\hat \bmu_{k-1}.$ 
By definition,
\[
  \hat \gamma_k ({\bf y}) = \min_{j \in {\mathcal P}_k ({\bf y})}   R_{k,j} ({\bf y}), 
\]
where
 \[
{\mathcal P}_k({\bf y}) = \{j  \mbox{ admissible and }
R_{k,j}({\bf y}) > 0\}.
\]
For admissible $j$, 
$R_{k,j}({\bf y}_0) \neq 0,$ and near ${\bf y}_0$  the functions ${\bf y}
\rightarrow R_{k,j}({\bf y})$ are continuous and of fixed sign. Thus, near ${\bf y}_0$ the
set ${\mathcal P}_k ({\bf y})$ stays fixed at $ {\mathcal P}_k ({\bf y}_0)$ and
(\ref{eq:singlept}) implies that
\[
  R_{k,k+1}({\bf y}) < R_{k,j} ({\bf y}), \qquad 
j > k+1, j \in {\mathcal P}_k ({\bf y}).
\]
Consequently, for ${\bf y}$ near ${\bf y}_0$, only variable $k+1$ joins the active set, and
so ${\mathcal A}_{k+1}({\bf y}) = {\mathcal A}_k \cup \{ k+1 \}$, and
\begin{equation}
    \label{eq:gammaky}
\hat \gamma_k ({\bf y}) = R_{k,k+1}({\bf y})
  = \frac{({\bf x}_l - {\bf x}_{k+1})'({\bf y} - 
\hat \bmu_{k-1} ({\bf y}))}{({\bf x}_l - {\bf x}_{k+1})'{\bf u}_k}.
  \end{equation}
This representation shows that both $\hat \gamma_k({\bf y})$ and hence $\hat
\bmu_k({\bf y}) = \hat \bmu_{k-1}({\bf y}) + \hat \gamma_k({\bf y}) {\bf u}_k$ are continuous
(resp. linear) near ${\bf y}_0$.

To show that $\hat \gamma_k$ is locally Lipschitz at ${\bf y}$, we set
$\bdelta = {\bf w} -
{\bf x}_{k+1}$ and write, using notation from (\ref{eq:lipdef}), 
\[
  \Delta \hat \gamma_k = \frac{\bdelta'(\Delta {\bf y} - \Delta \hat
  \bmu_{k-1})}{\bdelta' {\bf u}_k}.
\]
As ${\bf y}$ varies, there is a finite list of vectors $({\bf x}_l,{\bf x}_{k+1},{\bf u}_k)$
that can occur in the denominator term $\bdelta'{\bf u}_k$, and since all
such terms are positive [as observed below~(\ref{eq:singlept})], they
have a uniform positive lower bound, $a_{\min}$ say. Since $\Vert 
\bdelta \Vert
\leq 2$ and $\hat \bmu_{k-1}$ is Lipschitz $(L_{k-1})$ by assumption,
we conclude that
\[
  \frac{| \Delta \hat \gamma_k | }{ \Vert \Delta {\bf y} \Vert} \leq 
  2 a_{\min}^{-1} (1 + L_{k-1})  =: L_k.
\]
\upqed\end{pf}

\section{Consequences of the positive cone condition.}
\label{sec:cons-posit-cone}
\begin{lemma}
  \label{lem:rankone}
Suppose that $|{\mathcal A}_+| = |{\mathcal A} | + 1$ and that
$X_{{\mathcal A}+} = [ X_{{\mathcal A}} \ \ {\bf x}_+ ]$ 
\textup{(}where ${\bf x}_+ = s_j {\bf x}_j$ for some
$j \notin {\mathcal A}$\textup{)}.  Let $P_{\mathcal A} = 
X_{\mathcal A} G_{\mathcal A}^{-1}
X_{\mathcal A}'$ denote projection on ${\rm span}(X_{\mathcal A}),$ so that
$a = {\bf x}_+' P_{\mathcal A} {\bf x}_+ < 1.$ The  $+$-component 
of $G_{{\mathcal
    A}+}^{-1} {\bf 1}_{{\mathcal A}+}$ is 
\begin{equation}
  \label{eq:update}
  ( G_{{\mathcal A}+}^{-1} {\bf 1}_{{\mathcal A}+} )_+ = (1-a)^{-1} 
\biggl( 1 - \frac{{\bf x}_+'
  {\bf u}_{\mathcal A}}{A_{\mathcal A}} \biggr).
\end{equation}
Consequently, under the positive cone condition \textup{(\ref{eq:pcc})},
\begin{equation}
  \label{eq:goodangle}
  {\bf x}_+' {\bf u}_{\mathcal A} < A_{\mathcal A}.
\end{equation}
\end{lemma}

\begin{pf}
  Write $G_{{\mathcal A}+}$ as a partitioned matrix
\[
    G_{{\mathcal A}+} = 
  \pmatrix{
    X'X & X'{\bf x}_+ \cr  {\bf x}_+'X  & {\bf x}_+'{\bf x}_+   
    }
    =
    \pmatrix{
      A & B \cr B' & D}.
\]
Applying the formula for the inverse of a partitioned matrix
[e.g., \citet{rao73}, page 33],
\[
  ( G_{{\mathcal A}+}^{-1} {\bf 1}_{{\mathcal A}+} )_+ = 
- E^{-1} F' {\bf 1} + E^{-1},
\]
where
\begin{eqnarray*}
  E & =& D - B'A^{-1} B = 1 - {\bf x}_+' P_{\mathcal A} {\bf x}_+,   \\
  F & =& A^{-1}  B = G_{\mathcal A}^{-1} X' {\bf x}_+, 
\end{eqnarray*}
from which (\ref{eq:update}) follows. The positive cone condition
implies that $G_{{\mathcal A}+}^{-1} {\bf 1}_{{\mathcal A}+} > 0$, and so
(\ref{eq:goodangle}) is immediate.
\end{pf}

\section{Global continuity and Lemma~\ref{prop:regularity}.}
\label{sec:glob-cont-prop}

We shall call ${\bf y}_0$ a multiple point at step $k$ if two or more
variables enter at the same time. Lemma~\ref{lem:singlepoint} shows
that such points form a set of measure zero, but they can and do cause
discontinuities in $\hat \bmu_{k+1}$ at ${\bf y}_0$ in general.
We will see, however, that the positive cone condition prevents such
discontinuities.

We confine our discussion to double points, hoping that these
arguments will be sufficient to establish the same pattern of behavior
at points of multiplicity 3 or higher. In addition, by
renumbering, we shall suppose that indices $k+1$ and $k+2$ are those
that are added at double point ${\bf y}_0$. Similarly, for convenience only,
we assume that ${\mathcal A}_k({\bf y})$ is constant near ${\bf y}_0$.
Our task then is to show that, 
for ${\bf y}$ near a double point ${\bf y}_0$,
both $\hat \bmu_k({\bf y})$ and $\hat \bmu_{k+1}({\bf y})$ are continuous and
uniformly locally Lipschitz.

\begin{lemma}
  \label{lem:pcca}
  Suppose that ${\mathcal A}_k({\bf y}) = {\mathcal A}_k$ is constant near 
${\bf y}_0$ and
  that ${\mathcal A}_{k+} ({\bf y}_0) = {\mathcal A}_k \cup \{ k+1, k+2 \}$.  Then
  for ${\bf y}$ near ${\bf y}_0$, ${\mathcal A}_{k+}({\bf y}) \setminus 
{\mathcal A}_k$ can only
  be one of three possibilities, namely $\{ k+1 \},\{ k+2 \}$ or $\{
  k+1, k+2 \}$. In all cases $\hat \bmu_k({\bf y}) = \hat \bmu_{k-1} ({\bf y}) +
  \hat \gamma_k({\bf y}) {\bf u}_k$ as usual, and both $\gamma_k({\bf y})$ and
  $\hat \bmu_k({\bf y})$ are continuous and locally Lipschitz.
\end{lemma}

\begin{pf}
  We use notation and tools from the proof of Lemma
  \ref{lem:singlepoint}. Since ${\bf y}_0$ is a double point and the
  positivity set ${\mathcal P}_k({\bf y}) = {\mathcal P}_k$ near ${\bf y}_0$, we have
\[
    0 < R_{k,k+1}({\bf y}_0) = R_{k,k+2}({\bf y}_0) < R_{k,j}({\bf y}_0) 
\qquad
    \mbox{for }  j \in {\mathcal P}_k \setminus \{k+1, k+2\}.
\]
Continuity of $R_{k,j}$ implies that near ${\bf y}_0$ we still have
  \[
    0 < R_{k,k+1}({\bf y}), R_{k,k+2}({\bf y}) < \min \bigl\{ R_{k,j}({\bf y});  
     j \in {\mathcal P}_k \setminus \{k+1, k+2\} \bigr\}.
  \]
Hence ${\mathcal A}_{k+} \setminus {\mathcal A}_k$ must equal $\{k+1\}$ or $\{
k+2 \}$ or $\{ k+1, k+2 \}$ according as~$R_{k,k+1}({\bf y})$ is 
less than, greater than or equal to  $R_{k,k+2}({\bf y}).$ The continuity of
\[
  \hat \gamma_k ({\bf y}) = \min \{ R_{k,k+1}({\bf y}), R_{k,k+2}({\bf y}) \} 
\]
is immediate, and the local Lipschitz property follows from the
arguments of Lemma~\ref{lem:singlepoint}.
\end{pf}

\begin{lemma}
  \label{lem:doublepoint}
Assume the conditions of Lemma~\textup{\ref{lem:pcca}} and in addition that the
positive cone condition \textup{(\ref{eq:pcc})} holds. Then $\hat \bmu_{k+1}({\bf y})$
is continuous and locally Lipschitz near ${\bf y}_0$.
\end{lemma}

\begin{pf}
  Since ${\bf y}_0$ is a double point,  property (\ref{eq:singlept})
  holds, but now with equality when $j = k+1$ or $k+2$ and strict
  inequality otherwise. In other words, there exists $\delta_0 > 0$
  for which
  \[
    \hat C_{k+1}({\bf y}_0) - \hat c_{k+1,j}({\bf y}_0)  
    \cases{
      =0, & \quad \mbox{if }$j = k+2$, \cr
      \geq \delta_0,  &\quad  \mbox{if }$j > k+2$.}
     \]

Consider a neighborhood $B({\bf y}_0)$ of ${\bf y}_0$ and let $N({\bf y}_0)$ be the set
of double points in $B({\bf y}_0)$, that is, those for which ${\mathcal A}_{k+1}({\bf y})
\setminus {\mathcal A}_k = \{ k+1, k+2 \}$. We establish the convention that at
such double points $\hat \bmu_{k+1} ({\bf y}) = \hat \bmu_k({\bf y})$; at other
points ${\bf y}$ in $B({\bf y}_0)$, $\hat \bmu_{k+1} ({\bf y})$ is defined by $\hat
\bmu_k({\bf y}) + \hat \gamma_{k+1}({\bf y}) {\bf u}_{k+1}$ as usual.

Now consider those ${\bf y}$ near ${\bf y}_0$ for which ${\mathcal A}_{k+1}
({\bf y})
\setminus {\mathcal A}_k = \{k+1 \}$, and so, from the previous lemma,
${\mathcal A}_{k+2}({\bf y}) \setminus {\mathcal A}_{k+1} = \{ k+2 \}.$
For such ${\bf y}$, continuity and the local Lipschitz property for $\hat
\bmu_k$ imply that
\[
    \hat C_{k+1}({\bf y}) - \hat c_{k+1,j}({\bf y})     
\cases{
      = O(\Vert {\bf y} - {\bf y}_0 \Vert), &\quad  
\mbox{if }  $j = k+2$, \cr
      > \delta_0/2, &\quad  \mbox{if }  $j > k+2$.}
\]
It is at this point that we use the positive cone condition (via Lemma
\ref{lem:rankone}) to guarantee that $A_{k+1} > a_{k+1,k+2}.$ Also, since ${\mathcal A}_{k+1}({\bf y})
\setminus {\mathcal A}_k = \{k+1 \}$, we have
\[
  \hat C_{k+1} ({\bf y}) > \hat c_{k+1,k+2} ({\bf y}).
\]
These two facts together show that $k+2 \in {\mathcal P}_{k+1}({\bf y})$ and hence
that
\[
  \hat \gamma_{k+1} ({\bf y}) = \frac{\hat C_{k+1} ({\bf y}) - \hat c_{k+1,k+2}
  ({\bf y})}{A_{k+1} - a_{k+1,k+2}} = O( \Vert {\bf y} - {\bf y}_0 \Vert)
\]
is continuous and locally Lipschitz. In particular, as ${\bf y}$ approaches
$N({\bf y}_0)$, we have $\hat \gamma_{k+1} ({\bf y}) \rightarrow 0$.
\end{pf}

\renewcommand{\therem}{\Alph{section}.\arabic{rem}}
\begin{rem} 
We say that a function $g\dvtx
\mathbb{R}^n \rightarrow \mathbb{R}$ is \textit{almost differentiable}
if it is absolutely continuous on almost all line segments parallel to
the coordinate axes, and its partial derivatives (which consequently
exist a.e.) are locally integrable.
This definition of almost differentiability appears
superficially to be weaker than that given by Stein, but it is in fact
precisely the property used in his proof. Furthermore, this definition
is equivalent to the standard definition of weak differentiability
used in analysis.
\end{rem}

\begin{pf*}{Proof of Lemma~\ref{prop:regularity}}
We have shown explicitly that $\hat \bmu_k ({\bf y})$ is continuous and
uniformly locally Lipschitz near\ single and double points. Similar
arguments extend the property to points of multiplicity 3 and
higher, and so all points ${\bf y}$ are covered. Finally, absolute
continuity of ${\bf y} \rightarrow \hat \bmu_k ({\bf y})$ on line segments is a
simple consequence of the uniform Lipschitz property, and so $\hat
\bmu_k$ is almost differentiable.
\rightqed\end{pf*}
\end{appendix}

\section*{Acknowledgments.}

The authors thank Jerome Friedman, Bogdan Popescu, Saharon Rosset 
and Ji Zhu for helpful discussions.

\printaddresses

\begin{thebibliography}{99}

\bibitem[\protect\citeauthoryear{Breiman, Friedman, Olshen and Stone}{1984}]{BFOS84}
{\sc Breiman, L., Friedman, J., Olshen, R.} and {\sc  Stone, C.} 
(1984). {\it Classification and Regression Trees}.
Wadsworth, Belmont, CA.
\MR{726392}

\bibitem[\protect\citeauthoryear{Efron}{1986}]{efro:1986}
{\sc Efron, B.}  (1986). How biased is the apparent
  error rate of a prediction rule? {\em J. Amer. Statist.
 Assoc.} {\bf 81} 461--470.
\MR{845884}

\bibitem[\protect\citeauthoryear{Efron and Tibshirani}{1997}]{ET97}
{\sc Efron, B.} and  {\sc Tibshirani, R.}  (1997).
Improvements on cross-validation: The $.632+$ bootstrap method. {\em J. Amer.
Statist. Assoc.} {\bf 92} 548--560.
\MR{1467848}

\bibitem[\protect\citeauthoryear{Freund and Schapire}{1997}]{FS95}
{\sc Freund, Y.} and  {\sc Schapire, R.}  (1997).
A decision-theoretic generalization of online learning and an application to
boosting. {\it J.  Comput.  System Sci.} {\bf 55} 119--139.
\MR{1473055}

\bibitem[\protect\citeauthoryear{Friedman}{2001}]{Fr99}
{\sc Friedman, J.}  (2001). Greedy function
approximation: A gradient boosting machine. {\em Ann. Statist.}
\textbf{29} 1189--1232.
\MR{1873328}

\bibitem[\protect\citeauthoryear{Friedman, Hastie and Tibshirani}{2000}]{FHT00}
{\sc Friedman, J., Hastie, T.} and {\sc Tibshirani, R.} (2000). 
Additive logistic regression: A statistical view of
boosting (with discussion). {\em Ann. Statist.} {\bf 28} 337--407.
\MR{1790002}

\bibitem[\protect\citeauthoryear{Golub and Van Loan}{1983}]{GVL83}
{\sc Golub, G.} and {\sc Van Loan, C.}  (1983).
{\em Matrix Computations}. Johns Hopkins Univ. Press, Baltimore, MD.
\MR{733103}

\bibitem[\protect\citeauthoryear{Hastie, Tibshirani and Friedman}{2001}]{ElemStatLearn}
{\sc Hastie, T., Tibshirani, R.} and  {\sc Friedman, J.}  (2001).
{\em The Elements of Statistical Learning\textup{:} Data
Mining\textup{,} Inference and Prediction}. Springer, New York.
\MR{1851606}

\bibitem[\protect\citeauthoryear{Lawson and Hanson}{1974}]{LH74}
{\sc Lawson, C.} and {\sc Hanson, R.}  (1974).
{\em Solving Least Squares Problems}. 
Prentice-Hall, Englewood Cliffs, NJ.
\MR{366019}

\bibitem[\protect\citeauthoryear{Mallows}{1973}]{Ma73}
{\sc Mallows, C.}  (1973). Some comments on $C_p$.
{\em Technometrics} {\bf 15} 661--675.

\bibitem[\protect\citeauthoryear{Meyer and Woodroofe}{2000}]{meyer00}
{\sc Meyer, M.} and {\sc  Woodroofe, M.} (2000).
On the degrees of freedom in shape-restricted regression. 
{\em Ann.  Statist.} {\bf 28} 1083--1104.
\MR{1810920}

\bibitem[\protect\citeauthoryear{Osborne, Presnell and Turlach}{2000a}]{osborne00}
{\sc Osborne, M., Presnell, B.} and {\sc Turlach, B.}  (2000a). 
A new approach to variable selection in least squares
problems. {\it IMA J. Numer. Anal.} {\bf 20} 389--403. 
\MR{1773265}

\bibitem[\protect\citeauthoryear{Osborne, Presnell and Turlach}{2000b}]{osborne00:_lasso_dual}
{\sc Osborne, M.~R., Presnell, B.} and {\sc Turlach, B.}  (2000b). 
On the LASSO and its dual. {\it J.~Comput.  Graph. Statist.} {\bf 9} 319--337.
\MR{1822089}

\bibitem[\protect\citeauthoryear{Rao}{1973}]{rao73}
{\sc Rao, C.~R.}  (1973). {\it Linear Statistical
Inference and Its Applications}, 2nd ed. Wiley, New York.
\MR{346957}

\bibitem[\protect\citeauthoryear{Stein}{1981}]{St81}
{\sc Stein, C.}  (1981). Estimation of the mean of
a multivariate normal distribution. {\it Ann. Statist.} {\bf 9} 1135--1151.
\MR{630098}

\bibitem[\protect\citeauthoryear{Tibshirani}{1996}]{Ti96}
{\sc Tibshirani, R.} (1996). Regression shrinkage
and selection via the lasso. {\it J. Roy. Statist. Soc. Ser. B.} {\bf
58} 267--288.
\MR{1379242}

\bibitem[\protect\citeauthoryear{Weisberg}{1980}]{W80}
{\sc Weisberg, S.}  (1980). {\em Applied Linear
Regression}. Wiley, New York.
\MR{591462}

\bibitem[\protect\citeauthoryear{Ye}{1998}]{ye93:_in}
{\sc Ye, J.}  (1998). On measuring and correcting
the effects of data mining and model selection. {\em J. Amer.
Statist. Assoc.} \textbf{93} 120--131.
\MR{1614596}
\end{thebibliography}
\end{document}